\documentclass[12pt,leqno]{article}

\usepackage{amsmath}
\usepackage{amsfonts}
\usepackage{amssymb}
\usepackage{amsthm}
\usepackage{latexsym}

\def\g{\mathfrak{g}}
\def\a{\mathfrak{a}}
\def\b{\mathfrak{b}}

\def\k{\mathfrak{k}}

\def\p{\mathfrak{p}}
\def\n{\mathfrak{n}}
\def\m{\mathfrak{m}}
\def\l{\mathfrak{l}}

\def\s{\mathfrak{s}}

\def\N{\mathbb{N}}

\def\R{\mathbb{R}}
\def\C{\mathbb{C}}
\def\ad{\rm{\, ad\, }}
\def\Ad{\rm{\, Ad \,}}
\def\Om{\Omega}
\def\ind{{\rm ind}}

\def\d{\delta}
\def\si{\sigma}
\def\v{\varphi}
\def\e{\varepsilon}
\def\r{\rho}

\def\O{\Omega}
\def\o{\omega}
\def\De{\Delta}

\def\ve{\varepsilon}

\def\E{\mathcal{E}}

\def\M{\mathcal{M}}

\def\T{\mathcal{T}}

\def\cd{\cdot}
\def\iy{\infty}

\def\ol#1{\overline{#1}}

\def\cds#1#2{#1,\cdots,#2}

\def\hb#1{\hbox{#1}}
\def\va#1{\vert #1\vert}

\def\no#1{\Vert #1\Vert }
\def\opno#1{\Vert #1\Vert_{\mathrm op} }
\def\pa#1{\{#1\}}

\def\res#1{_{\vert #1}}
\def\inv{^{-1}}

\def\hb #1{\hbox{#1}}

\def\iy{\infty}

\def\ol#1{\overline{#1}}

\def\cd{\cdot}
\def\cds{\cdots}

\def\hb#1{\hbox{#1}}
\def\va#1{\vert #1\vert}

\def\pa#1{\{#1\}}

\def\ca#1{\mathcal{#1}}

\def\Log#1{\rm{Log}({\it #1})}

\def\Im{\mathrm{\, Im \,}}
\def\Re{\mathrm{\, Re\, }}
\def\stacksx{\tilde\xi}

\def\xiphi{\xi\otimes^p_{\rho}\phi}
\def\etapsi{\eta\otimes^{p'}_{\rho}\psi}

\def\l({\left(}
\def\r){\right)}
\def\lan{\langle}
\def\ran{\rangle}

\def\sp#1#2{\langle #1,#2\rangle }

\newtheorem{thm}{Theorem}[section]

\newtheorem{exs}{Examples}[section]
\newtheorem*{rem}{Remark}
\newtheorem*{rems}{Remarks}

\newtheorem*{cor1}{Corollary of Theorem \ref{intermediatetheorem}}
\newtheorem*{cor2}{Corollary of Proposition \ref{limite}}
\newtheorem{lemme}{Lemma}[section]
\newtheorem{proposition}{Proposition}[section]

\begin{document}

\title{Holomorphic \(L^{p}\)-type for sub-Laplacians on connected Lie groups}
\author{Jean Ludwig, Detlef M\"uller, Sofiane Souaifi
\thanks{
This work has been supported by the the IHP network HARP ``Harmonic
Analysis and Related Problems'' of the European Union\newline
{\em keywords:} connected Lie group, sub-Laplacian, functional
calculus, $L^p$-spectral multiplier,  symmetry,
holomorphic $L^p$-type, $L^p$-transference by induced representations\newline
{\em 2001 Mathematics Subject Classification:} 22E30, 22E27, 43A20
}
}
\date{}

\maketitle

\begin{abstract}
We study the problem of determining all connected Lie groups
$G$ which have the following property (hlp): every sub-Laplacian
$L$ on $G$ is of holomorphic $L^p$-type for $1\leq p<\infty, \ p\ne 2.$
First we show that semi-simple non-compact Lie groups with finite
center have this property, by using holomorphic families of
representations in the class one principal series of $G$ and the
Kunze-Stein phenomenon. We then apply  an  $L^p$-transference
principle, essentially due to Anker,  to show that every connected Lie
group $G$ whose  semi-simple quotient by its radical is non-compact  has
property (hlp). For the convenience of the reader, we give a
self-contained  proof of this  transference principle, which
generalizes the well-known Coifman-Weiss principle.   One is thus
reduced to  studying those groups for which the semi-simple quotient is
compact, i.e.~to compact extensions of solvable Lie groups. In this
article, we consider  semi-direct extensions of exponential solvable Lie
groups by connected compact Lie groups. It had been proved  in
\cite{hlm} that every exponential solvable Lie group
$S$, which has a non-$*$ regular co-adjoint orbit whose restriction to the
nilradical is closed, has  property (hlp), and we show here that
(hlp)  remains valid for compact extensions of these groups.
\hspace{-0.2cm}
\end{abstract}

\tableofcontents

\section{Introduction}\label{introduction}
\setcounter{equation}{0}
A comprehensive discussion of the problem studied in this article,
background information and references to further literature can be
found in \cite{l-m}. We shall therefore content ourselves in this
introduction by recalling some notation and results from
\cite{l-m}.\\
Let $(X,d\mu)$ be a measure space. If $T$ is
a self-adjoint linear operator on a $L^2$-space
$L^2(X,d\mu)$, with spectral resolution $T=\int_{\R}
\lambda dE_\lambda$, and if $F$ is a bounded Borel function on
$\R$, then we call $F$ an $L^p$-{\em multiplier for } $T$ 
($1\le p<\infty$), if $F(T):=\int_\R F(\lambda)dE_\lambda$ extends from
$L^p\cap L^2(X,d\mu)$ to a bounded operator on $L^p(X,d\mu)$. We
shall denote by ${\cal M}_p(T)$ the space of all $L^p$-multipliers
for $T$, and by $\sigma_p(T)$ the $L^p$- spectrum of $T.$ We say
that $T$ is of {\em holomorphic $L^p$-type}, if there exist some
non-isolated point $\lambda_0$ in the $L^2$-spectrum $\sigma_2(T)$
and an open complex neighborhood ${\cal U}$ of $\lambda_0$ in
$\C$, such that every $F\in {\cal M}_p(T)\cap C_\infty (\R)$
extends holomorphically to ${\cal U}$. Here, $C_\infty (\R)$
denotes the space of all continuous functions on $\R$ vanishing at
infinity.\\
Assume in addition that there exists a linear subspace ${\cal D}$
of $L^2(X)$ which is $T$-invariant and dense in $L^p(X)$ for every
$p\in [1,\infty[$, and that $T$ coincides with the closure of its
restriction to ${\cal D}$. Then, if $T$ is of holomorphic
$L^p$-type, the set ${\cal U}$ belongs to the $L^p$-spectrum of
$T$, i.e.
\begin{equation*}
\overline{\cal U}\subset\sigma_p(T).
\end{equation*}
In particular,
\begin{equation*}
\sigma_2(T)\subsetneq \sigma_p(T).
\end{equation*}
Throughout this article, $G$ will denote a connected Lie group.\\
Let $dg$ be a {\em left-invariant} Haar measure on $G$. If
$(\pi,\mathcal{H}_\pi)$ is a unitary representation of $G$ on
the Hilbert space ${\cal H}={\cal H}_\pi$, then we denote the
integrated representation of $L^1(G)=L^1(G,dg)$ again by $\pi$,
i.e.~$\pi(f)\xi:=\int _G f(g)\pi(g)\xi\,dg$ for every $f\in
L^1(G), \ \xi\in{\cal H}.$ For $X\in \g$, we denote by $d\pi(X)$
the infinitesimal generator of the one-parameter group of unitary
operators $t\mapsto \pi(\exp tX)$. Moreover, we shall often identify $X$
with the corresponding {\em right-invariant} vector field
$X^r f(g):= \lim_{t\to 0} \frac 1 t [f((\exp tX)g)-f(g)]$ on $G$ and
write $X=X^r.$\\

Let $X_1,\,\ldots,\,X_k$ be elements of $\g$ which generate $\g$ as a Lie
algebra, which just means that the corresponding right-invariant vector
fields satisfy  H\"ormander's condition. The corresponding sub-Laplacian
$L:=-\sum_{j=1}^k X_j^2$ is then essentially self-adjoint on 
${\cal D}(G)\subset L^2(G)$ and hypoelliptic.  Denote by
$\{e^{-tL}\}_{t>0}$ the heat semigroup generated by $L$. Since $L$
is right $G$-invariant, for every $t>0$, $e^{-tL}$ admits a
convolution kernel $h_{t}$ such that
\begin{displaymath}
e^{-tL}f= h_{t}\ast f,
\end{displaymath}
where $\ast$ denotes the usual convolution product in $L^{1}(G)$. The
function $(t,g)\mapsto h_{t}(g)$ is smooth on $\R_{>0}\times G,$
since the differential operator $\frac{\partial}{\partial t}+L$ is
hypoelliptic. Moreover, by \cite[Theorem VIII.4.3 and
Theorem V.4.2]{vsc}, the heat kernel $h_t$ as well as its
right-invariant derivatives admit Gaussian type estimates in terms
of the Carnot-Carath\'eodory distance $\delta$ associated to the
H\"ormander system $X_1,\dots X_k.$ \\
In particular,  for every right-invariant differential
operator $D$ on $G$, there exist constants $c_{D,t}, C_{D,t}>0$,
such that
\begin{equation}\label{gauss-estimate}
\vert D h_t(g)\vert \leq C_{D,t} e^{-c_{D,t}
\delta(g,e)^{2}},\quad \textrm{ for all }g\in G, t>0.
\end{equation}
Let now $F_0\in{\cal {M}}_p(L)$. By duality, we may assume that 
$1\le p\le 2.$ With $F_0$, also the function
$\lambda\mapsto F(\lambda):=e^{-\lambda}F_0(\lambda)$ lies in 
${\cal M}_p(L)$, since
$F(L)=e^{-L}F_0(L),$ where the heat operator $e^{-L}$ is
bounded on every $L^p(G)\, (1\leq p<\infty).$ Now by \cite[Lemma 6.1]{l-m},
the operator $F_0(L)$ is bounded also on all the spaces
$L^q(G),\, p\leq q\leq p'$. Hence for every test function $f$ on
$G$,
\begin{displaymath}
\begin{array}{c}
F(L)(f)=F_0(L)(e^{-L}(f))=F_0(L)(h_1*f)\\
=F_0(L)(h_{1/2}*h_{1/2}*f)=(F_0(L)h_{1/2})*h_{1/2}*f,
\end{array}
\end{displaymath}
by the right invariance of the operator $F_0(L)$. Since $h_{1/2}$ is
contained in every $L^q(G),\ 1\leq q\leq \infty$, in particular in
$L^1(G)$, we see that the operator $F(L)$ acts by convolution from
the left with the function $(F_0(L)h_{1/2})*h_{1/2}$ which is
contained in every $L^q(G),\, p\leq q\leq p',$ and so are all its
derivates from the right. We can thus identify the operator $F(L)$
with the $C^\infty$-function
$F(L)\delta:=(F_0(L)h_{1/2})*h_{1/2},$ i.e.
\begin{equation*}
F(L)(f)=(F(L)\delta)*f, \quad f\in \bigcup_{p\le q\le p'} L^q(G).
\end{equation*}
%
\noindent Recall that the {\em modular function} $\Delta_G$ on $G$
is defined by the equation
$$
\int_G f(xg)dx=\Delta_G (g)^{-1} \int_Gf(x)dx, \qquad g\in G.
$$
We put:
\begin{eqnarray*}
\check{f}(g) & := & f(g^{-1}),\nonumber\\
f^*(g) & := & \Delta_G^{-1}(g) \overline {{f(g^{-1})}}.
\end{eqnarray*}
Then $f\mapsto f^*$ is an isometric involution on $L^1(G)$, and
for any unitary representation $\pi$ of $G$, we have:
\begin{equation}\label{adjoint}
\pi(f)^*=\pi(f^*)\ .
\end{equation}
The group $G$ is said to be {\em symmetric}, if the associated
group algebra $L^1(G)$ is symmetric, i.e.~if every element $f\in
L^1(G)$ with $f^*=f$ has a real spectrum with respect to the
involutive Banach algebra $L^1(G).$\\
In this paper we consider connected Lie groups for which
every sub-Laplacian is of holomorphic $L^p$-type. First, in the
Section \ref{sectbig2},  we consider connected semi-simple Lie groups $G$
with finite center. We construct a holomorphic family of
representations $\pi_{(z)}$ of $G$ on mixed $L^p$-spaces (see
Section \ref{sect3}). Applying these representations to $h_1$, we
obtain a holomorphic family of compact operators on these spaces
(see Section \ref{sect4}). Using the Kunze-Stein phenomenon on
semi-simple Lie groups (see Section \ref{sect5}), the eigenvectors
of the operators $\pi_{(z)}(h_1)$ allow us to construct a holomorphic
family of $L^p$-functions on $G$ which are eigenvectors for $F(L)$, if
$F\in {\cal M}_p(T)\cap C_\infty (\R).$ From the corresponding
 holomorphic family of eigenvalues we can read off that $F$ admits a
holomorphic extension in a neighborhood of some element in the
spectrum of $L$ (see Section \ref{sect6}). This gives us:
\begin{thm}\label{maintheorem1}
Let $G$ be a non compact connected semi-simple
Lie group with finite center. Then every sub-Laplacian on $G$ is
of holomorphic $L^p$-type, for $1\le p<\infty, \ p\ne 2.$
\end{thm}
\begin{rem}
{\rm
Even if at the end of the proof, we consider only ordinary $L^p$-spaces, 
we need representations on mixed $L^p$-spaces. They are used to get some 
isometry property and then to apply the Kunze-Stein phenomenon.
}
\end{rem}
In Sections \ref{induced} and \ref{transfer}, 
we discuss respectively p-induced representations and a
generalization of the Coifman-Weiss transference principle
\cite{coifman-weiss}. We consider
a separable locally compact group
$G$, and  an isometric representation
$\rho$  of a closed subgroup $S$ of $G$ on  spaces of $L^p$-type, 
e.g.~$L^p$-spaces $L^p(\Omega).$ Denote by $\pi_p:=\ind_{p,S}^G\,\rho$
the $p$-induced representation of $\rho.$  We prove, among other
results,  that, for any function $f\in L^1(G),$ the operator norm
of $\pi_p(f)$ is bounded by the norm of the convolution operator
$\lambda_G(f)$ on $L^p(G),$ provided the group $S$ is amenable.
Here, $\lambda_G$ denotes the left-regular representation. It
should be noted that we do not require the group $G$ to be
amenable.  As an   application
we obtain the $L^p$-transference of a convolution operator on $G$
to a convolution operator on the quotient group $G/S,$ in the case
where $S$ is an amenable closed, normal  subgroup.\\
When preparing this article, we were not aware of J.-Ph.~Anker's article
\cite{anker} which, to a large extent, contains these transference results,
and which we also recommend for further references to this topic. We are
indebted to N. Lohou\'e for
informing us on  Anker's work \cite{anker} as well as on the influence of
C. Herz on the
development of this field (compare \cite{herz}).   For the
convenience of the reader, we have nevertheless decided to include our
approach to
these transference results, since it  differs from Anker's by the use of
a suitable cross section for $G/S,$  which we feel makes the arguments a
bit easier.\\

Applying this transference principle, we obtain the following
generalization of Theorem \ref{maintheorem1} in Section \ref{sectbig4}:
\begin{thm}\label{maintheorem2} Let $G=\exp\g$ be a connected
Lie group, and denote by  $S=\exp{\s}$  its radical. If $G/S$ is not
compact, then every sub-Laplacian    on $G$ is of
holomorphic $L^p$-type, for any $1\leq p<\infty, \ p\ne 2$.
\end {thm}
It then suffices to study connected Lie groups for which $G/S$ is
compact. 
In Section \ref{sectbig5}, we shall consider groups $G$ which are the  semi-direct product of a compact
group $K$ with a non-symmetric exponential solvable group $S$ from a certain class.
The exponential solvable non-symmetric Lie groups have been completely
classified by Poguntke \cite{p} (with previous contributions by Leptin,
Ludwig and Boidol) in terms of a purely Lie-algebraic condition
\eqref{boidolcond}. Let us describe this condition, which had been first
introduced by Boidol in a different context \cite{boidol}.\\
Recall that the unitary dual of $S$ is in one to one
correspondence with the space of coadjoint orbits in the dual space $\s^*$
of $\s$ via the Kirillov map, which associates with a given point
$\ell\in\s^*$ an irreducible unitary representation $\pi_{\ell}$ (see,
e.g., \cite[Section 1]{hlm}).\\
If $\ell$ is an element of $\s^*$, denote by
$$\s(\ell):=\{X\in \s|\,\ell([X,Y])=0,\textrm{ for all }Y\in \s\}$$
the {\em stabilizer } of $\ell$ under the coadjoint action
$\ad^*$. Moreover, if $\m$ is any Lie algebra, denote by
$$\m=\m^1\supset \m^2\supset \dots$$
the descending central series of $\m$, i.e.~$\m^2=[\m,\m]$, and
$\m^{k+1}=[\m,\m^k]$. Put
$$\m^\infty =\bigcap_k\m^k.$$
Then $\m^\infty$ is the smallest ideal $\k$ in $\m$ such that $\m/\k$
is nilpotent. Put
$$\m(\ell):=\s(\ell)+[\s,\s].$$
Then we say that $\ell$ respectively the associated coadjoint
orbit $\Om(\ell):=\Ad^*(G)\ell$ satisfies {\em Boidol's
condition}\, \eqref{boidolcond}, if
\begin{equation}\label{boidolcond}
\ell\mid_{\m(\ell)^\infty} \ne 0.\tag{B}
\end{equation}
According to \cite{p}, the group $S$ is non-symmetric if and only
if there exists a coadjoint orbit satisfying Boidol's condition.\\
If $\Om$ is a coadjoint orbit, and if $\n$ is the nilradical of
$\s$, then
$$\Om|_{\n}:=\{\ell|_{\n}: \ell\in\Om\}\subset \n^*$$
will denote the restriction of $\Om$ to $\n.$\\
We show that the methods developped in \cite{hlm} can also be
applied to the case of a compact extension of an exponential
solvable group and thus obtain
\begin{thm}\label{maintheorem3}
Let $G=K\ltimes S$ be a semi-direct product of a compact Lie group
$K$ with an exponential solvable Lie group $S$, and assume that
there exists a coadjoint orbit $\Om(\ell)\subset \s^*$ satisfying
Boidol's condition, whose restriction to the nilradical $\n$ is
closed in $\n^*$. Then every sub-Laplacian on $G$ is of
holomorphic $L^p$-type, for $1\le p<\infty,\, p\ne 2.$
\end{thm}
%
\begin{rems}${}$
{\rm
\begin{itemize}
\item[{\bf (a)}] A sub-Laplacian $L$ on $G$ is of
holomorphic $L^p$-type if and only if every continuous bounded
multiplier $F\in{\cal M}_p(L)$ extends holomorphically to an open
neighborhood of a non-isolated point in $\sigma_2(L).$
\item[{\bf (b)}] If the restriction of a coadjoint orbit to the nilradical is
closed, then the orbit itself is closed (see \cite[Thm. 2.2]{hlm}).
\item[{\bf (c)}] What we really use in the proof is the following property of
the orbit $\Om:$
\begin{quote}
\noindent {\em $\Om$ is closed, and for every real character $\nu$
of $\s$ which does not vanish on $\s(\ell),$ there exists a
sequence $\{\tau_n\}_n$ of real numbers such that
$\lim_{n\to\infty}(\Om+\tau_n\nu)=\infty$ in the orbit space.}
\end{quote}
This property is a consequence of the closedness of $\Om|_\n.$
There are, however, many examples where the condition above is
satisfied, so that the conclusion of the theorem still holds, even
though the restriction of $\Om$ to the nilradical is not closed
(see e.g.~\cite[Section 7]{hlm}). We do not know whether the condition above
automatically holds whenever the orbit $\Om$ is closed.\\
Observe that, contrary to the semisimple case, we need to consider 
representations on mixed $L^p$-spaces till the end of the proof.
\end{itemize}
}
\end{rems}
In all the sequel, if $M$ is a topological space, $C_0(M)$ will mean 
the space of compactly supported continuous funcions on $M$.\\ 
As usual, if $S$ is a Lie group, $\s$ will denote its Lie algebra.


\section{The semi-simple case}\label{sectbig2}
\setcounter{equation}{0}

\subsection{Preliminaries}\label{sect2}
If $E$ is a vector space, denote by $E^*$ its algebraic dual. If
it is real, $E_{\C}$ denotes its complexification. Let $F$ be a
vector subspace of $E$. We identify in the sequel the restriction
$\lambda_{| F}$ of $\lambda\in E^*$ or
$E_{\C}^*$ to an element of respectively $F^*$ or $F_{\C}^*$. \\
Let $G$ be a connected semisimple real Lie group with finite
center and $\g$
its Lie algebra.\\
Fix a Cartan involution $\theta$ of $G$ and denote by $K$ the
fixed point group for $\theta$. The Cartan decomposition of the
Lie algebra $\g$ of $G$ with respect to $\theta$ is given by
$$\g=\k\oplus \p,$$
where $\k$ is the Lie algebra of $K$ and $\p$ the $-1$-eigenspace
in $\g$ for
the differential of $\theta$, denoted again by $\theta$.\\
We fix a subspace $\a$ of $\p$ which is maximal with respect to
the condition that it is an abelian subalgebra of $\g$. It is
endowed with the scalar product $(\cdot,\cdot)$ given by the
Killing form $B$, which is positive definite on $\p$. By duality,
we endow $\a^*$ with the corresponding, induced scalar product,
which we also denote by $(\cdot,\cdot)$. Let $|\cdot|$ be the
associated norm on $\p$ and
$\a^*$.\\
For any root $\alpha\in\a^*, $ we denote by $\g_{\alpha}$ the
corresponding root space, 
i.e.~$\g_{\alpha} :=\{X\in\g|\,[H,X]=\alpha(X)\textrm{ for }H\in \a\}.$ 
We fix a set $R^+$
of positive roots of $\a$ in $\g$. Let $P$ denote the
corresponding minimal parabolic subgroup of $G$, containing
$A:=\exp \a$, and $P=MAN$ its
Langlands decomposition.\\
Denote by $\rho$ the linear form on $\a$ given by
$$\rho(X):= \frac{1}{2}\, \mathrm{tr}\,(\mathrm{ad}X_{|\n})\textrm{ for all }
X\in \a,$$
where $\n$ is the Lie algebra of $N$.\\
%
Let $\Vert \cdot \Vert$ denote the ``norm'' on $G$ defined in
\cite[\S 2]{vdbs}. Recall that, for $g\in G$, $\Vert g\Vert$ is
the operator norm of $\mathrm{Ad}\, g$ considered as an operator
on $\g$, endowed with the real Hilbert structure, 
$(X,Y)\mapsto -B(X,\theta Y)$ as scalar product. This norm is $K$-biinvariant
and, according to \cite[Lemma 2.1]{vdbs}, satisfies the following
properties:
\begin{equation}\label{norme}
{\setlength\arraycolsep{0pt}
\begin{array}{l}
\Vert \cdot\Vert \textrm{ is a continuous and proper function on }G,\\
\Vert g\Vert = \Vert\theta(g) \Vert = \Vert g^{-1}\Vert \geq 1;\\
\Vert xy\Vert \leq \Vert x \Vert \, \Vert y\Vert;\\
\textrm{there exists } c_1, c_2 >0\textrm{ such that, for }Y\in
\p\textrm{, then }
e^{c_1 |Y|} \leq \Vert \exp Y\Vert \leq e^{c_2 |Y|};\\
\textrm{for all }\,\,\,a\in A, n\in N,\,\Vert a\Vert\leq \Vert an\Vert.
\end{array}
}
\end{equation}
%
We choose a basis for $\a^*$, following, for example, 
\cite[p.~220]{cowling}.\\
Let $\alpha_{1},\dots,\alpha_{r}$ denote the simple roots in
$R^{+}$. By the Gram-Schmidt process, one constructs from the
basis $\{\alpha_{1},\dots,\alpha_{r}\}$ of $\a^{*}$ an orthonormal
basis $\{\beta_{1},\dots,\beta_{r}\}$ of $\a^{*}$ in a such way
that, for every $j=1,\dots,r$, the vector space
$\mathrm{Vect}\{\beta_{1},\dots,\beta_{j}\}$ spanned by
$\{\beta_{1},\dots,\beta_{j}\}$ agrees with
$\mathrm{Vect}\{\alpha_{1},\dots,\alpha_{j}\}$, and, for every
$1\leq k< j\leq r$, $(\beta_{j},\alpha_{k})=0$. Define $H_{j}$
($j=1,\dots,r$) as the element of $\a$ given by
$\beta_{k}(H_{j})=\delta_{jk}$ ($k=1,\dots,r$), and put
\begin{displaymath}
\begin{array}{l}
\a_{j} := \R H_{j}\\
\a^{j} := \sum_{k=1}^{j}\a_{k}\\
R^{j} := R^{+}\cap \mathrm{Vect}\{\alpha_{1},\dots,\alpha_{j}\}\\
R_{j} := R^{j}\setminus R^{j-1} \textrm{, with } R^{0} := \emptyset\\
\n^{j} := \sum_{\alpha\in R^{j}}\g_{\alpha}\\
\n_{j} := \sum_{\alpha\in R_{j}}\g_{\alpha}.
\end{array}
\end{displaymath}
%
We define, for $j=1,\dots,r$, the reductive Lie subalgebra
$\m^{j}$ of $\g$ by setting
\begin{displaymath}
\m^{j} := \theta(\n^{j}) + \m + \a^{j} + \n^{j}.
\end{displaymath}
In this way, we obtain a
finite sequence of reductive Lie subalgebras of $\g,$
\begin{displaymath}
\m=: \m^{0} \subset \m^{1} \subset \dots \subset \m^{r} = \g ,
\end{displaymath}
such that
\begin{displaymath}
\m^{j} = \theta(\n_{j}) + \m^{j-1} + \a_{j} + \n_{j} \,\quad (j=1,\dots,r).
\end{displaymath}
Then $\m^{j-1}\oplus \a_{j}\oplus \n_{j}$ is a parabolic
subalgebra
of $\m^{j}$.\\
%
Observe that $G$ is a real reductive Lie group in the
Harish-Chandra class (see e.g.~\cite[p.~58]{g-v}, for the
definition of this class of reductive Lie groups). We can then
inductively define a decreasing sequence of reductive real Lie
groups $M^{j}$ in the Harish-Chandra class, starting from 
$M^r = G$, in the following way.\\ Let $P_j$ denote the parabolic
subgroup of $M^{j}$ corresponding to the parabolic subalgebra
$\m^{j-1}\oplus\a_{j}\oplus \n_{j}$, and $P_j = M^{j-1}A_j N_j$
its Langlands decomposition. Here $A_j$ (resp.~$N_j$) is the
analytic subgroup of $M^j$ with Lie algebra $\a_j$ (resp.~$\n_j$),
$M^{j-1}A_j$ is the centralizer in $M^j$ of $A_j$, and
$$M^{j-1}:=\bigcap_{\chi\in \mathrm{Hom}(M^{j-1}A_j,\R_{+}^{\times})}\mathrm{Ker}\, \chi$$
(see e.g.~\cite[Theorem 2.3.1]{g-v}). \\
Moreover, $M^{j-1}A_j$ normalizes $N_j$ and $\theta(N_j),$ and
$M^{j-1}$ is a reductive Lie subgroup of $M^j$, in the
Harish-Chandra class, with Lie algebra $\m^{j-1}$ (see \cite[Proposition 2.1.5]{g-v}).\\
%
Put $K^{j} := M^{j}\cap K$ ($j=1,\dots,r$). Then $K^{j}$ is the
maximal compact subgroup of $M^{j}$ related to the Cartan involution
$\theta_{\vert M^{j}}$ of $M^{j}$ (see e.g.~\cite[Theorem
2.3.2, p.~68]{g-v}). Hence, $M^{j}$ is the product
\begin{displaymath}
K^j P_j = K^{j}M^{j-1}A_{j}N_{j}.
\end{displaymath}
%
Fix invariant measures $dk$, $dm$, $da$, $dn$, $dm_{j}$, $dk_{j}$,
$da_{j}$, $dn_{j}$ for respectively $K$, $M$, $A$, $N$, $M^{j}$,
$K^{j}$,
$A_{j}$, $N_{j}$.\\
Choose an invariant measure $dx$ on $G$ such that
\begin{equation}\label{mesure}
\int_G \varphi(x)\, dx= \int_{K\times A\times N}
a^{2\rho}\varphi(kan)\,dkdadn,\ \textrm{ for all}\ \varphi\in
C_0(G)
\end{equation}
(see, e.g., \cite[Proposition 2.4.2]{g-v}).\\
%
We shall next recall an integral formula. Let $S$ be a reductive
Lie group in the Harish-Chandra class, and let $S=K\exp \p$ be its
Cartan decomposition, where $K$ is a maximal compact subgroup of
$S$. Let $Q$ be a parabolic subgroup of $S$ related to the above
Cartan decomposition, and let $Q=M_Q A_Q N_Q$ be its Langlands
decomposition.\\ 
Let $K_Q:=K\cap Q = K\cap M_Q,$ and put, for $k\in K$, $[ k ]:= k K_Q$
in $K/K_Q$. We extend this notation to $S$ by putting, for
$s=kman$, $(k,m,a,n)\in K\times M_Q\times A_Q\times N_Q$ , $[s]:=kK_Q$. This is still well defined even though the
representation of $s$ in $KM_Q A_Q N_Q$ is not unique. In fact,
$$s =kman =k'm'a'n'$$
if and only if $a'=a$, $n'=n$, and $k'=kk_Q$, $m'=k_Q ^{-1}m$ for
some $k_Q\in K_Q$ (see e.g.~\cite[Theorem 2.3.3]{g-v}). From this
we see that the decomposition above becomes unique, if we require
$m$ to be in
$\exp(\m_Q\cap \p)$.\\
Every $s\in S$ thus admits a unique decomposition $s=kman,$ with
$(k,m,a,n)\in K\times \exp(\m_Q\cap \p)\times A_Q\times N_Q$. We
then write $k_Q(s):=k$, $m_Q(s):=m$, $a_Q(s):=a$ and $n_Q(s):=n,$
i.e.
$$s=k_Q(s)m_Q(s)a_Q(s)n_Q(s).$$
In particular, $[s]=k_Q(s)K_Q$.\\ For $y\in S$ and $k\in K$, we
define $y[k]\in K/K_Q$ as follows:
$$y[k]:=[yk].$$
Moreover, for any $\gamma\in\a^*_{\C}$ and $Y\in\a,$ we put 
$(\exp Y)^\gamma:=e^{\gamma(Y)}.$
\begin{lemme}\label{change}
Fix an invariant measure $dk$ on $K$ and let $d[k]$ denote the
corresponding left invariant measure on $K/K_Q$. For any $y\in S$,
we then have
$$d(y[k]) = a_Q(yk)^{-2\rho_Q} d[k],$$
where $\rho_Q\in \a_Q^*$ is given by
$\rho_Q(X)=\frac{1}{2}\mathrm{tr}\, (\mathrm{ad} X_{|\n_Q})$
($X\in \a_Q$); that is, for any $f\in C(K/K_Q)$,
$$\int_{K/K_Q} f([k])d[k] = \int_{K/K_Q} a_Q(yk)^{-2\rho_Q}f(y[k])d[k].$$
\end{lemme}
\begin{proof}
We follow the proof of \cite[Proposition 2.5.4]{g-v}. Let 
$f\in C(K/K_Q)$. Consider $f$ as a right $K_Q$-invariant function on
$K$. Choose $\chi\in C_0(M_Q A_Q N_Q)$ such that
$$\int_{M_Q\times A_Q\times N_Q} a^{2\rho_Q}\chi(man)dmdadn = 1$$
and $\chi(k_Q q) =\chi(q),$ for all $k_Q\in K_Q,\, q\in Q,$ and
put, for $s\in S$ with $s=kman$, $(k,m,a,n)\in K\times M_Q\times A_Q\times N_Q$,
$$h(s):=f(k)\chi(man).$$
Notice that the function $h$ is well defined, independently of the
chosen decomposition $s=kman$ of $s$, since
$f$ (respectively $\chi$) is right (respectively left) $K_Q$-invariant.\\
Let $dm$, $da$ and $dn$ be invariant measures on $M_Q$, $A_Q$ and
$N_Q$, respectively. Choose an invariant measure $dx$ on $S$ such
that
$$\int_S \varphi(x)\, dx= \int_{K\times M_Q\times A_Q\times N_Q}
a^{2\rho_Q}\varphi(kman)\,dkdmdadn,\ \textrm{ for all }\varphi\in C_0(S)$$ 
(see, e.g., \cite[Proposition 2.4.3]{g-v}).
Then
$$\int_{K}f(k)dk=\int_S h(x)\, dx.$$
On the other hand, by left invariance,
$$\int_S h(x)\, dx = \int_S h(yx)\, dx=\int_{K\times M_Q\times A_Q\times N_Q} a^{2\rho_Q}h(ykman)\,dkdmdadn.$$ 
Write $yk = k_Q(yk)m_Q(yk)a_Q(yk)n_Q(yk)$. Since the elements of $M_Q$ and
$A_Q$ commute, we get
$$ykman = k_Q(yk)m_Q(yk)ma_Q(yk)an_Q(yk)^{(ma)^{-1}}n,$$
where $n_Q(yk)^{ma}=ma n_Q(yk)(ma)^{-1}\in N_Q,$ since $M_Q A_Q$
normalizes $N_Q.$ Therefore,
$$h(ykman) = f([yk])\chi(m_Q(yk)ma_Q(yk)an_Q(yk)^{(ma)^{-1}}n).$$
We thus obtain, by left invariance of $dm$, $da$, $dn$,
\begin{displaymath}
\begin{array}{rl}
\int_S h(yx)\,dx & = \int_{K\times M_Q\times A_Q\times N_Q}
a^{2\rho_Q}a_Q(yk)^{-2\rho_Q}f([yk])\chi(man)\,dkdmdadn\\
& = \int_K a_Q(yk)^{-2\rho_Q}f(y[k])\,dk.
\end{array}
\end{displaymath}
The lemma follows.
\end{proof}
\begin{rem}
{\rm
Let $P=MAN$ be a minimal parabolic subgroup of $G$ contained in $Q$.\\
If we decompose $s\in S$ via the Iwasawa decomposition $S=KAN$ as
$$s=k(s)a(s)n(s),$$
where $k(s)\in K$, $a(s)\in A$ and $n(s)\in N$, we can check that
$k(s) = k_Q(s)$ and $a(s) = a(m_Q(s))\,a_Q(s)$, where $a(m_Q(s))$
lies in fact in $\exp(\m_Q \cap \a)$. Since this space is
orthogonal to $\a_Q$ with respect to the scalar product
$(\cdot,\cdot)$ on $\p$, for any $\lambda\in \a_Q^*$ we have
$\lambda|_{\m_Q\cap\a}=0,$ hence
$$a(s)^{\lambda} = a_Q(s)^{\lambda}.$$
With these considerations, the lemma above can also be deduced,
for example, from \cite[Lemma 2.4.1]{wallach}.
}
\end{rem}
%

We return now to our semisimple Lie group $G.$ In the sequel, we
shall use another basis of $\a^*,$ given as follows. For
$j=1,\dots,r$, let $\rho_{j}$ denote the element of $\a_{j}^*$
defined by
$$\rho_{j}(X):= \frac{1}{2} \mathrm{tr}\,(\mathrm{ad}X_{|\n_{j}})\textrm{ for all }X\in \a_{j}.$$ 
Notice that we can identify $\rho_{j}$ with
the restriction $\rho_{|\a_{j}}$ of
$\rho$ to $\a_{j}$.\\
By \cite[ Lemma 4.1]{cowling}, $\rho_{j}$ and $\beta_{j}$ are
scalar multiples of each other. In particular, the family
$\{\rho_{j}\}$ is an orthogonal basis of $\a^*$, and therefore of
$\a_{\C}^*$. For every $\nu\in\a^*$ (resp.~$\nu\in \a_{\C}^*$),
define $\nu_{j}$ ($j=1,\dots,r$) in $\R$ (resp.~$\C$) by the
following:
$$\nu=\sum_{j=1}^r \nu_{j}\rho_{j}.$$
Recall that, for $j=1,\dots,r$, $P_j$ is a parabolic subgroup of
the reductive real Lie group $M^j$, which lies in the
Harish-Chandra class. Put therefore, by taking
$(S,K,Q):=(M^j,K^j,P_j)$ in the discussion above, $k_j :=k_{P_j}$,
$m_{j-1} := m_{P_j}$, $a_j := a_{P_j}$ and $n_j := n_{P_j}$.
Then, any $g\in M^j$ has a unique decomposition
$g=k_j(g)m_{j-1}(g)a_j(g)n_j(g)$, with $k_j(g)\in K^j$,
$m_{j-1}(g)\in \exp(\m^{j-1}\cap \p)$, $a_j(g)\in A_j$ and
$n_j(g)\in N_j$. Notice that $\m_0\cap\p=\{0\}$, i.e.~$m_0(g)=e.$
%
\begin{lemme}\label{borne}
Denote by $r_{y}$ the right multiplication with $y\in G$.
Let $j\in \{1,\dots,r\}$, $g\in M^{j}$ and $k_{l}\in K^{l}$ 
($l=1,\dots,j$).\\
We define recursively the element $g_l$ of $M^l$, $l=1,\dots,j$,
starting from $l=j$, by putting $g_j:=g$ and 
$g_{l-1} := m_{l-1}(g_lk_l),$ i.e.
$$g_l=m_l\circ(r_{k_{l+1}}\circ m_{l+1})\circ \cdots \circ (r_{k_j}\circ
m_j)(g),\quad 1\le l\le j-1.$$ 
Then, the following estimate holds:
\begin{equation*}
\Vert \Pi_{l=j}^1 a_l(g_lk^l)\Vert \leq \Vert g\Vert.
\end{equation*}
\end{lemme}
\begin{proof}
We first show that, for $1\leq p\leq j$,
\begin{equation}\label{eq1.1}
\Vert g\Vert=\Vert \Pi_{l=j}^{p} a_l(g_lk_l)\cdot
m_{p-1}(g_pk_p)\cdot \Pi_{l=p}^{j}n_l (g_lk_l)k_l^{-1}\Vert .
\end{equation}
(Here the products are non-commutative products, in which the
order of the factors is indicated by the order of indices.) We use
an induction, starting from $p=j$. If $p=j$ and $g\in M^j,$ then
$$\Vert g\Vert=\Vert g k_jk_j^{-1}\Vert= \Vert
k_j(gk_j)m_{j-1}(gk_j)a_j (gk_j)n_j(gk_j)k_j^{-1}\Vert.$$ 
Using the left $K$-invariance of the norm and the fact that
$a_j(gk_j)\in A_j$ and $m_{j-1}(gk_j)\in M^{j-1}$ commute, we find
that
$$\Vert g\Vert=\Vert a_j(gk_j)m_{j-1}(gk_j)n_j(gk_j)k_j^{-1}\Vert,$$ 
so that \eqref{eq1.1} holds for $p=j.$ Assume now, by induction, that
\eqref{eq1.1} is true for $p+1$ in place of $p,$ i.e.
$$\Vert g\Vert=\Vert \Pi_{l=j}^{p+1} a_l(g_lk_l)\cdot
m_{p}(g_{p+1}k_{p+1})\cdot \Pi_{l=p+1}^{j}n_l
(g_lk_l)k_l^{-1}\Vert .$$ 
We then decompose
$$m_{p}(g_{p+1}k_{p+1})k_p=g_pk_p=k_p(g_pk_p)m_{p-1}
(g_pk_p)a_p(g_pk_p)n_p(g_pk_p).$$ Since 
$k_p(g_pk_p)\in K^p\subset M^l,$ for $p\le l\le j,$ 
it commutes with $a_l(g_lk_l), $ for
$l=p+1,\dots, j,$ and therefore, because of the $K$-invariance of
$\Vert\cdot\Vert,$ we have
$$ \Vert g\Vert=\Vert \Pi_{l=j}^{p+1} a_l(g_lk_l)\cdot
m_{p-1}(g_{p}k_{p}) a_p(g_pk_p)n_p(g_pk_p) k_p^{-1}\cdot
\Pi_{l=p+1}^{j}n_l (g_lk_l)k_l^{-1}\Vert .$$ Moreover,
$a_p(g_pk_p)$ commutes with $m_{p-1}(g_pk_p),$ and so
\eqref{eq1.1}
follows.\\
Applying now \eqref{eq1.1} for $p=1$, we obtain
\begin{equation}\label{eq1.2}
\Vert \Pi_{l=j}^{1}
a_l(g_lk_l)\Pi_{l=1}^{j}n_l(g_lk_l)k_l^{-1}\Vert = \Vert g\Vert.
\end{equation}
By right $K$-invariance of the norm, the left hand side of this
equation is equal to
$$\Vert \Pi_{l=j}^{1} a_l(g_lk_l)\Pi_{l=1}^{j}n_l(g_lk_l)k_l^{-1}\Pi_{l'=j}^{1}k_{l'}\Vert.$$ 
Notice that we can write
$\Pi_{l=1}^{j}n_l(g_lk_l)k_l^{-1} \Pi_{l'=j}^{1}k_{l'}$ as
follows:
$$n_1(g_1k_1)\big(k_1^{-1}n_2(g_2k_2)k_1\big)\big((k_2k_1)^{-1}
n_3(g_3k_3)k_2k_1\big)\cdots
\big((\Pi_{l=j-1}^{1}k_l)^{-1}n_j(g_jk_j)\Pi_{l'=j-1}^{1}k_{l'}\big).$$
Since $(\Pi_{l'=p-1}^{1}k_{l'})^{-1}$, $2\leq p\leq j$, lies in
$K^{p-1}\subset M^{p-1}$ and thus normalizes $N_p$, we get that
$$\Pi_{l=1}^{j}n_l(g_lk_l)k_l^{-1} \Pi_{l'=j}^{1}k_{l'}\in N.$$
Using the last property of the norm given in \eqref{norme}, the
left hand side of \eqref{eq1.2} is then greater or equal to 
$\Vert\Pi_{l=j}^{1} a_l(g_lk_l)\Vert,$ which proves the lemma.
\end{proof}

\subsection{A holomorphic family of representations of $G$ on mixed
$L^{p}-$ spaces}\label{sect3}
For $\nu\in \a_{\C}^{*}$, let $\mathcal{M}(G,P,\nu)$ denote the
space of complex valued measurable functions $f$ on $G$ satisfying
the following covariance property:
\begin{displaymath}
f(gman) = a^{-(\nu+\rho)}f(g)\textrm{ for all }g\in G,\,m\in
M,\,a\in A, \,n\in N.
\end{displaymath}
The space $\mathcal{M}(G,P,\nu)$ is endowed with the left regular
action of $G$, denoted by $\tilde{\pi}_{\nu}$, i.e.,
$[\tilde{\pi}_\nu(g)f](g')=f(g^{-1}g')$. The representations
$\tilde{\pi}_\nu$ form the class-one principal series.\\ Let
$\mathcal{M}(K/M)$ be the space of right $M$-invariant measurable
functions on $K$.\\ The restriction to $K$ of functions on $G$
gives us a linear isomorphism from $\mathcal{M}(G,P,\nu)$ onto
$\mathcal{M}(K/M)$. Denote by $I_\nu : f\mapsto f_{\nu}$ the
inverse mapping. Then $f_{\nu}\in \mathcal{M}(G,P,\nu)$ is given
by
\begin{displaymath}
f_{\nu}(kan) := a^{-(\nu+\rho)}f(k)\textrm{ for all }k\in K,\,
a\in A,\,n\in N,
\end{displaymath}
if $G=KAN$ denotes the Iwasawa decomposition of $G$.\\
If we intertwine the representation $\tilde{\pi}_\nu$ with
$I_\nu$, we obtain a representation $\pi_\nu$ of $G$ on
$\mathcal{M}(K/M)$, given by
\begin{displaymath}
(\pi_{\nu}(g)f)_{\nu} = \tilde{\pi}_{\nu}(g)f_{\nu}\textrm{, if }f\in \mathcal {M}(K/M),\, g\in G.
\end{displaymath}
Denote by $d\dot{k}_{j}$, for $j=1,\dots,r$, the quotient measure
on $K^{j}/K^{j-1}$ coming from $dk_{j}$. It is invariant by left
translations.
Notice that $K^{j-1}=K^j\cap M^{j-1}.$\\
We choose a left invariant measure $d\dot{k}$ on $K/M$ such that,
for any $f\in C(K/M),$
\begin{equation}\label{mesurek}
\int_{K/M}f(k)d\dot{k} = \int_{K^r/K^{r-1}}\cdots\int_{K^1/M}
f(k_r\cdots k_1)d\dot{k}_1\cdots d\dot{k}_r.
\end{equation}
Let $\underline{p}=(p_{1},\dots,p_{r})\in [1,+\infty [^{r}$.\\
One can easily see that, for every $f\in \mathcal{M}(K/M)$, 
$k'\in K$, the function on $K^j$ given by
\begin{displaymath}
k\mapsto \Big(\int_{K^{j-1}/K^{j-2}}\cdots \big(\int_{K^{1}/M}
\big| f(k'kk_{j-1} \cdots k_{1})\big|^{p_{1}}
d\dot{k}_{1}\big)^{p_{2}/p_{1}}\cdots d\dot{k}_{j-1}
\Big)^{1/p_{j-1}},
\end{displaymath}
is right $K^{j-1}$-invariant.\\
We can thus define the mixed $L^p$-space $L^{\underline{p}}(K/M)$,
as the space of all (equivalent classes of) functions $f$ in
$\mathcal{M}(K/M)$ whose mixed $L^p$-norm
\begin{displaymath}
\Vert f \Vert_{\underline{p}}:= \Big(\int_{K^{r}/K^{r-1}}\cdots
\big(\int_{K^{1}/M} \big| f(k_{r}\cdots k_{1})\big|^{p_{1}}
d\dot{k}_{1}\big)^{p_{2}/p_{1}}\cdots d\dot{k}_{r} \Big)^{1/p_{r}}
\end{displaymath}
is finite, endowed with this norm. This definition extends to the
case where
some of the $p_j$ are infinite, by the usual modifications.\\
Let $d$ denote the right $G$- and left $K$-invariant metric on
$G,$ given by
\begin{displaymath}
d(g,g') := \frac{1}{c_{1}}\log \Vert g'g^{-1} \Vert \,\,\,\,\, (g,\, g' \in G),
\end{displaymath}
where $c_1$ is the positive constant appearing in \eqref{norme}.
Notice that $d(g,e)=0$ if and only if $g$ lies in the center of
$G.$ In particular, $d$ is not separating.\\
Then, for $a={\rm exp} Y,$ with $Y\in\a\subset\p,$ and
$\gamma\in\a^*,$ we have, in view of the fourth property of
$\Vert\cdot\Vert$ in \eqref{norme}, that
\begin{equation}\label{2.5}
a^\gamma=|e^{\gamma(Y)}|\leq e^{|\gamma|\,|Y|}\leq
(e^{c_1|Y|})^{|\gamma|/c_1}\leq \Vert
a\Vert^{\frac{|\gamma|}{c_1}} =e^{|\gamma|\, d(a,e)}.
\end{equation}
\begin{proposition}\label{prop1}
For every $f\in L^{\underline{p}}(K/M)$ and $g\in G$, we have
\begin{displaymath}
\Vert \pi_{\nu}(g)f\Vert_{\underline{p}} \leq e^{|\sum_{j} (\frac{2}{p_j}-\mathrm{Re}\,\nu_{j}- 1)\rho_j| d(g,e)} 
\Vert f\Vert_{\underline{p}}.
\end{displaymath}
Thus $\pi_{\nu}$ defines a representation
$\pi_{\nu}^{\underline{p}}$ of $G$ on $L^{\underline{p}}(K/M)$.
Furthermore, this gives us an analytic family
$\{\pi_{\nu}^{\underline{p}}\}_{\nu\in \a_{\C}^{*}}$ of
representations of $G$ on $L^{\underline{p}}(K/M)$.
\end{proposition}
Before giving the proof, we show the following statement. We keep
the same notations as in Lemma \ref{borne}.
\begin{lemme}\label{lemint}
Let $g\in M^j,k\in K$ and $f_\nu\in\mathcal M(G,P,\nu).$ Then
\begin{displaymath}
{\setlength\arraycolsep{0pt}
\begin{array}{rcl}
& \Big(\int_{K^{j}/K^{j-1}}\cdots \big(\int_{K^{1}/M} \big|
f_{\nu}(kgk_{j} \cdots k_{1})\big|^{p_{1}}
d\dot{k}_{1}\big)^{p_{2}/p_{1}}\cdots d\dot{k}_{j} \Big)^{1/p_{j}} &\\
= &\Big(\int_{K^{j}/K^{j-1}}\cdots\big(\int_{K^{1}/M} \big|
\Pi_{l=j}^{1}a_l(g_lk_l)^{-(\mathrm{Re}\,\nu_{l}+1) \rho_{l}}
f_{\nu}(k\Pi_{l=j}^{1}k_{l}(g_lk_l)& )\big|^{p_{1}}
d\dot{k}_{1}\big)^{p_{2}/p_{1}}\cdots d\dot{k}_{j}
\Big)^{1/p_{j}}.
\end{array}
}
\end{displaymath}
\end{lemme}
\begin{proof}
We use induction on $j$. For $j=0$, one has, by right
$M$-invariance of $f$ and since $g\in M^0 = M$, that
\begin{displaymath}
|f_{\nu}(k g)|=|f_{\nu}(k)|.
\end{displaymath}
Assume that the statement is true for $j-1$. Observe that $a_j(gk_j)$ commutes with $k_{j-1}\cdots k_{1}\in M^{j-1},$ and
that 
$(k_{j-1} \cdots k_{1})^{-1} n_{j}(gk_j)k_{j-1} \cdots k_{1}\in N.$ 
Therefore, the covariance property of $f_\nu$ 
applied to the integration over $K^j/K^{j-1}$, implies
\begin{displaymath}
{\setlength\arraycolsep{0pt}
\begin{array}{rcl}
& \Big(\int_{K^{j}/K^{j-1}}\cdots \big(\int_{K^{1}/M} \big| f_{\nu}(kgk_{j} \cdots k_{1})\big|^{p_{1}}
d\dot{k}_{1}\big)^{p_{2}/p_{1}}\cdots d\dot{k}_{j} \Big)^{1/p_{j}} &\\
= & \Big(\int_{K^{j}/K^{j-1}} \cdots \big(\int_{K^{1}/M}\big| a_j(gk_j)^{-(\mathrm{Re}\,\nu_{j}+1)\rho_{j}}
f_{\nu}(kk_{j}(gk_j)m_{j-1}(gk_j) & k_{j-1} \cdots k_{1})\big|^{p_{1}} \\ & & d\dot{k}_{1}\big)^{p_{2}/p_{1}}\cdots
d\dot{k}_{j} \Big)^{1/p_{j}}.
\end{array}
}
\end{displaymath}
The statement holds, using the induction hypothesis, since $g=g_j$
and $m_{j-1}(gk_j)=g_{j-1}\in M^{j-1}$.
\end{proof}
\begin{proof}[Proof of Proposition \ref{prop1}]
If we apply \eqref{2.5} to 
$\gamma:=\sum_{j=1}^r (\frac 2{p_j}-\mathrm{Re}\,\nu_{j}- 1)\rho_j$ 
and notice that the
$\a_l$'s are pairwise orthogonal with respect to $(\cdot,\cdot)$,
we get, in view of Lemma \ref{borne}:
\begin{displaymath}
\sup_{k_j\in K^j,\,j=1,\dots,r} 
\Pi_{j=1}^{r}a_j(g_jk_j)^{(\frac{2}{p_j}-\mathrm{Re}\,\nu_{j}- 1)\rho_j}\leq \Vert g\Vert^{\frac{|\gamma|}{c_1}}=e^{|\gamma|d(g,e)}.
\end{displaymath}
On the other hand, according to the above lemma and Lemma
\ref{change}, applied successively to the integrations over
$K^j/K^{j-1}$, $j=1,\dots,r$, we have
\begin{displaymath}
\Vert \pi_{\nu}(g^{-1})f\Vert_{\underline{p}}\leq
\sup_{k_j\in K^j,\,j=1, \dots,r}
(\Pi_{j=r}^{1}a_j(g_jk_j)^{(\frac{2}{p_j}-\mathrm{Re}\,\nu_{j}-1)\rho_j} \big)\Vert f\Vert_{\underline{p}}.
\end{displaymath}
The first assertion of the proposition is now evident, since
$d(g^{-1},e)=d(g,e).$\\
In order to prove the analyticity of the family of representations
$\pi_{\nu}^{\underline{p}},$ choose
$\underline{p}=(p_1,\dots,p_r)\in [1, \infty[^r$ and denote by
$\underline{p}'=(p_1',\dots,p_r')\in ]1,\infty]^r$ the tuple of
conjugate exponents, i.e., $1/p_j +1/p_j'=1$. Then, for $f\in L^{\underline{p}}(K/M)$ , 
$u\in L^{\underline{p}'}(K/M)=(L^{\underline{p}}(K/M))'$ and $g\in G$,
we have
\begin{displaymath}
\langle\pi_{\nu}^{\underline{p}}(g)f,u\rangle=\int_{K/M}
(\pi_{\nu}^{\underline{p}} (g)f)(k)\overline{u(k)}d\dot{k}=
\int_{K/M}a(g^{-1}k)^{-(\nu+\rho)}f(\kappa(g^{-1}k))
\overline{u(k)}d\dot{k}.
\end{displaymath}
Here, the functions $ \kappa(\cdot), a(\cdot), n(\cdot)$ on $G$
are given by the unique factorization $g=\kappa(g)a(g)n(g)$ of
$g,$
according to the Iwasawa decomposition $G=KA N$.\\
Obviously, the expression above is analytic in $\nu\in\a^*_{\C},$
which finishes the proof.
\end{proof}
For $t=(t_1,\dots,t_r)\in ]0,+\infty[^{r},$ let
$\Omega_{t}:=\{\nu\in \a_{\C}^{*}\big|\,|\mathrm{Re}\,\nu_{j} |< t_{j}\ \textrm{for all}\ j=1,\dots,r\}$. 
Moreover, for $p\ge 0$,
let $\overline p:=(p,\dots,p)\in \R^r.$
\begin{proposition}\label{proprep} ${}$
\begin{itemize}
\item[{\bf (i)}] For all $\underline{p}\in [1,+\infty [^{r}, \,f\in L^{\underline{p}} (K/M),\, \nu\in \Omega_{t},\, g\in G,$ we
have
$$\Vert \pi_{\nu}^{\underline{p}}(g) f\Vert_{\underline{p}} \leq e^{\sum_{j}(t_{j} + 1)|\rho_j| \, d(g,e)} 
\Vert f\Vert_{\underline{p}.}
$$
\item[{\bf (ii)}] Let $\nu \in \a_{\C}^{*}$, and let
$\underline{q}$ be an element of $[1,+\infty[^{r}$ satisfying
\begin{displaymath}
\mathrm{Re}\, \nu_{j} = \frac 2{q_{j}} - 1,\quad j=1,\dots,r.
\end{displaymath}
Then, for all $g\in G, \, f\in L^{\underline{q}}(K/M),$
\begin{displaymath}
\Vert \pi_{\nu}^{\underline{q}}(g) f\Vert_{\underline{q}} = \Vert f\Vert_{\underline{q}}.
\end{displaymath}
Furthermore, for $\nu\in i\a^{*}$, $\pi_{\nu}^{\overline{2}}$ is a
unitary representation of $G$.
\end{itemize}
\end{proposition}
\begin{proof}
(i) results immediately from the estimate given in Proposition
\ref{prop1} and (ii) from Lemmas \ref{lemint} and \ref{change},
since, for such $\underline{q}$, we have
$a^{-(\mathrm{Re}\,\nu_l+1)\rho_l}= a^{-2\rho_l/q_l},$ if $a\in A_l$.
\end{proof}

\subsection{A holomorphic family of compact operators}\label{sect4}
Let $L=-\sum_1^k X_j^2$ be a fixed sub-Laplacian on $G$.
The estimate \eqref{gauss-estimate}, in combination with the
estimate in Proposition \ref{proprep}\,(i), easily implies that
the operator
$$\pi_{\nu}^{\underline{p}}(h_1)f :=\int_G h_1(x)\pi_{\nu}^{\underline{p}}(x)f\,dx,\,\, f\in L^{\underline{p}}(K/M),$$ 
is well defined and bounded on
$L^{\underline{p}}(K/M)$. In fact these operators are even
compact. To see this, let us put, for $\nu\in \Omega_{1}$ and
$k_1,\,k_2\in K,$
\begin{eqnarray}\label{kernel}
K_{\nu}(k_1,k_2):=c_{G}\int_{M\times A\times N}
a^{-\nu+\rho}h_{1}(k_1 (man)^{-1}k_{2}^{-1})dmdadn,
\end{eqnarray}
where $c_{G}$ is the positive constant given by $d(x^{-1})= c_{G}dx$ (which
exists, since $G$ is unimodular).\\
\begin{lemme}\label{Kcontinu}
The integral in \eqref{kernel} is absolutely convergent and
defines a continuous, right $M$-invariant kernel function on
$K\times K,$ i.e.~$K_\nu(k_1m',k_2 m')=K_\nu(k_1,k_2)$ for every
$m'\in M.$
\end{lemme}
\begin{proof}
In order to prove that the integral in \eqref{kernel} is
absolutely convergent, we put
$$I:=\int_{M\times A\times N}
|a^{-\nu+\rho}h_{1}(k_1 (man)^{-1}k_{2}^{-1})|dmdadn.
$$
Then, in view of \eqref{gauss-estimate}, we have
$$I \leq C\int_{M\times A\times
N}a^{-\mathrm{Re}\,\nu+\rho}e^{-cd(k_1(man)^{-1}k_{2}^{-1},e)^2}dmdadn.$$
Using the $K$-bi-invariance of the norm $\Vert\cdot \Vert$ on $G$
and the inclusion $M\subset K$, we get that
$$ d(k_1 (man)^{-1}k_{2}^{-1},e) = d(kan,e)\textrm{, for any }k\in K.$$
Moreover, by \eqref{2.5} and \eqref{norme},
$$a^{-2\rho} a^{-\mathrm{Re}\,\nu+\rho}=a^{-\mathrm{Re}\,\nu-\rho}\leq
e^{|\mathrm{Re}\,\nu +\rho|d(kan,e)}.$$
We thus get, since
$|\mathrm{Re}\,\nu +\rho|\leq 2\sum_j |\rho_j |$ for $\nu\in \Omega_1$,
$$I \leq C \int_{K\times A\times N} a^{2\rho}e^{2\sum_j |\rho_j
|d(kan,e)}e^{-cd(kan,e)^2}dkdadn,$$ for every $k_1,\,k_2\in K$,
which is in fact equal to
$$C \int_G e^{2\sum_j |\rho_j |d(x,e)}e^{-cd(x,e)^2}\,dx.$$ Since
$G$ is unimodular and has exponential volume growth, it is easy to
see that this integral is finite. Moreover, since the integrand in
\eqref{kernel} depends continuously on $k_1$ and $k_2,$ we see
that $K_\nu$ is continous.\\
In order to prove the right $M$-invariance of $K_\nu,$ let $m'\in M$. One has, for any 
$(m,a,n)\in M\times A\times N$,
$$(man)^{m'} = m^{m'}a n^{m'}.$$
According to the invariance of $dm$, we then have, for any
$k_1,\,k_2\in K$,
$$K_{\nu}(k_1m',k_2m') = c_{G}\int_{M\times A\times N}
a^{-\nu+\rho}h_{1}(k_1 (man^{m'})^{-1}k_{2}^{-1})dmdadn.$$
Furthermore, it is easy to check that, for any $\varphi\in C_0(N)$,
$$\int_N \varphi(n^{m'})\,dn = \int_N \varphi(n)\,dn.$$ Indeed, since
$G=KAN$, there exists $\phi\in C_0(G)$
such that
$$\varphi(n) = \int_{K\times A} a^{2\rho}\phi(kan)\,dkda.$$
According to our choice of the Haar measure $dx$ on $G$
(c.f.~\eqref{mesure}), we have
$$\int_G\phi(x)\,dx = \int_{K\times A\times N}a^{2\rho}\phi(kan)\,dkdadn
=\int_N \varphi(n)\, dn,$$ 
and using the invariance of $dx$ and
$dk$, in combination with the commutation and normalization
properties of $m'\in M$, we see that
$$
\int_N \varphi(n)\,dn =\int_G \phi(x^{m'})\, dx= \int_{K\times A\times N}a^{2\rho} \phi(kan^{m'})\,dkdadn 
= \int_N\varphi(n^{m'})\,dn.
$$
We thus conclude that
$$K_{\nu}(k_1m',k_2m')= K_{\nu}(k_1,k_2).$$

\end{proof}
Put, for $\nu\in\Omega_{1}$,
\begin{displaymath}
T(\nu) := \pi_{\nu}(h_{1}).
\end{displaymath}

\begin{proposition}\label{compactop}
The operator $T(\nu)$ is represented by the integral kernel
$K_\nu,$ i.e.
\begin{eqnarray}\label{intker}
(T(\nu)f)(k_1 M)=\int_{K/M} K_\nu(k_1,k) f(k)\, d\dot{k}, \ f\in
L^1(K/M).
\end{eqnarray}
In particular, $T(\nu)=\pi_{\nu}^{\underline{p}}(h_1)$ is a
compact operator on every mixed $L^p$-space
$L^{\underline{p}}(K/M),$ $\underline{p}\in [1,+\infty[^{r},$
which we then shall also denote by $T_{\underline{p}}(\nu),$ in
order to indicate the space on which it acts. Moreover, the family
of compact operators $\nu\mapsto T_{\underline{p}}(\nu)$ is
analytic (in the sense of Kato
\cite{kato}) on $\Omega_{1}.$\\
Furthermore, for $\nu\in i\a^{*}$, the operator
$T_{\overline{2}}(\nu)$ is self-adjoint on $L^{\overline{2}}(K/M)$.
\end{proposition}
\begin{proof}
We have, by definition, for any $k_1\in K$, that
$$(T(\nu)f)(k_1) = \int_G h_1(x)(\pi_{\nu}(x)f)(k_1)\,dx.$$
By invariance of $dx$, this is equal to
$$c_G \int_G h_1(k_1 x^{-1})f_{\nu}(x)\,dx.$$
Now, according to our choice of $dx$ and using the covariance
property of $f_\nu$, we obtain
$$(T(\nu)f)(k_1) = c_G \int_{K\times A\times N}a^{2\rho}a^{-(\nu+\rho)} h_1(k_1(an)^{-1}k^{-1})f_\nu(k)\,dkdadn.$$ 
Since
$dk$ is invariant and $M\subset K$, this can be written, using
also the right $M$-invariance of $f_\nu$, as follows:
$$(T(\nu)f)(k_1) = c_G \int_{K\times M\times A\times N}a^{-\nu+\rho} h_1(k_1(man)^{-1}k^{-1})f_\nu(k)\,dkdmdadn.$$
\\
But, $f_\nu = f$ on $K,$ and thus \eqref{intker} follows, by
Fubini's
theorem.\\
Since the kernel $K_\nu$ is continuous on the compact space
$K\times K,$ by Lemma \ref{Kcontinu}, it follows that
$T_{\underline{p}}(\nu)$ is a compact operator on
$L^{\underline{p}}(K/M),$ and the analytic dependence of $K_\nu$,
which is evident from \eqref{intker}, implies that, for every
$\underline{p}\in [1,+\infty[^{r},$ the family of operators
$T_{\underline{p}}(\nu)$ is analytic on $\Omega_1.$\\
Finally, if $\nu\in i\a^*,$ then $\pi_\nu^{\overline{2}}$ is
unitary, and since $h_1(x)=\overline{h_1(x^{-1})}, $ we see (from 
\eqref{adjoint}) that
the operator $\pi_\nu^{\overline{2}}(h_1)$ is self-adjoint.
\end{proof}

\subsection{Some consequences of the Kunze-Stein phenomenon}\label{sect5}
Observe that, by H\"older's inequality, for any $\underline{p}\in [1,2]^r$ and any $\underline{q}\in [\underline{p},\underline{p}']$, 
we have
\begin{equation}\label{eq9.0}
\Vert f \Vert_{\underline{q}}\leq \Vert f \Vert_{\underline{p}'},
\quad\textrm{ for all }f\in L^{\underline{p}'}(K/M),
\end{equation}
since the compact space $K/M$ has normalized measure $1.$
Therefore,
$L^{\underline{p}'}(K/M)$ is a subspace of $L^{\underline{q}}(K/M)$.\\
Notice also that $L^p(K/M)= L^{\overline{p}}(K/M)$ and
$\Vert\cdot\Vert_{\overline{p}}=\Vert\cdot\Vert_p, $ by our choice
of measure on $K/M$ (c.f.~\eqref{mesurek}).\\
As a consequence of the Kunze-Stein phenomenon (see \cite{k-s} and
\cite{cowling}), we shall prove:
\begin{proposition}\label{interpol}
Let $1< p_0 < 2$ and $\nu_0\in \a^*\setminus\{0\}$. There exist
$\varepsilon >0$ and $C> 0$, such that, for any $\xi,\,\eta\in L^{{p}_0}(K/M)$ and $z\in \C$ with $|\mathrm{Re}\,z|<\varepsilon$,
\begin{equation}\label{eq9}
\Vert\langle\pi_{z\nu_0}(\cdot)\xi,\eta\rangle\Vert_{L^{p_0'}(G)}
\leq C\Vert \xi\Vert_{{p}_0'}\,\Vert \eta\Vert_{{p}_0'}.
\end{equation}
\end{proposition}
\begin{proof}
Observe that, for every $\nu\in i\a^*$, the representation
$\pi_\nu$ is unitarily equivalent to $\tilde{\pi}_\nu$. Therefore,
given $\delta>0$, we obtain from \cite{cowling}, that there is a
constant $C_\delta >0$, such that, for any $2+\delta \leq r'\leq\infty$ and $\xi,\, \eta\in L^{{2}}(K/M)$, we have:
\begin{equation}\label{eq10}
\Vert \langle\pi_\nu (\cdot)\xi,\eta\rangle\Vert_{L^{r'}(G)}\leq
C_\delta \Vert\xi\Vert_{{2}}\Vert \eta\Vert_{{2}}\textrm{,\quad
provided }\mathrm{Re}\,\nu =0.
\end{equation}
Indeed, in \cite{cowling}, this is only stated for $\nu =0$, but
the proof
easily extends to arbitrary $\nu\in i\a^*$.\\
On the other hand, we have the estimate:
\begin{equation}\label{eq11}
\Vert\langle \pi_\nu (\cdot)\xi,\eta\rangle\Vert_{L^\infty(G)}\leq
\Vert\xi\Vert_{\underline{q}}\Vert\eta\Vert_{\underline{q}'},\quad
\underline{q}\in [1,+\infty[^r,
\end{equation}
for any $\xi\in L^{\underline{q}}(K/M)$, $\eta\in L^{\underline{q}'}(K/M)$, provided that:
\begin{equation}\label{eq12}
\mathrm{Re}\,\nu_j = \frac{2}{q_j} - 1,\,\, j=1,\dots,r.
\end{equation}
This is an immediate consequence of Proposition
\ref{proprep}\,(ii), since $\pi_{\nu}^{\underline{q}}$ is
isometric, if \eqref{eq12} is satisfied.\\ Let $\theta_0\in ]0,1[$
be given by $\frac{2}{p_0} = 1+\theta_0$. Since ${2}\in [{p}_0,{p}_0']$ and since $\underline{q}\in [\overline{p}_0,\overline{p}_0']$, 
if $\underline{q}$ satisfies
\eqref{eq12} and $|\mathrm{Re}\, \nu_j |\leq \theta_0$,
$j=1,\dots,r$, we can unify \eqref{eq10} and \eqref{eq11}, using
\eqref{eq9.0}, as follows.\\ Given $\delta >0$, there is a
constant $C_\delta \geq 1$ such that, for all $\xi\in L^{{p}_0}(K/M)$, $\eta\in L^{{p}_0'}(K/M)$:
\begin{displaymath}
\begin{array}{c}
\Vert \langle \pi_\nu (\cdot)\xi,\eta\rangle\Vert_{L^{r'}(G)}\leq
C_\delta \Vert\xi\Vert_{{p}_0'}\Vert
\eta\Vert_{{p}_0'}\textrm{,\quad for all }r'\in
[2+\delta,+\infty]\textrm{, if }\mathrm{Re}\,\nu =0,
\end{array}
\end{displaymath}
and
\begin{displaymath}
\begin{array}{c}
\Vert\langle \pi_\nu (\cdot)\xi,\eta\rangle\Vert_{L^\infty(G)}\leq
\Vert\xi\Vert_{{p}_0'}\Vert\eta\Vert_{{p}_0'}\textrm{,\quad if }|\mathrm{Re}\, \nu_j|\leq \theta_0,\, j=1,\dots, r.
\end{array}
\end{displaymath}
If we choose $\nu=z\nu_0$, and put $\Psi_z :=\langle
\pi_{z\nu_0}(\cdot)\xi,\eta\rangle$, for $\xi,\,\eta\in
L^{{p}_0'}(K/M)$ fixed, this implies that
\begin{displaymath}
\Vert\Psi_{iy}\Vert_{L^{r'}(G)}\leq C_\delta
\Vert\xi\Vert_{{p}_0'}\Vert \eta\Vert_{{p}_0'}\textrm{,\quad for all }r'\in [2+\delta,+\infty]\textrm{ and }y\in \R,
\end{displaymath}
and
\begin{displaymath}
\Vert\Psi_{\pm\theta_1+iy}\Vert_{L^{\infty}(G)}\leq C_\delta
\Vert\xi\Vert_{{p}_0'}\Vert \eta\Vert_{{p}_0'}\textrm{,\quad for all }y\in \R,
\end{displaymath}
if we put $\theta_1 := \theta_0 /\max_{j=1,\dots,r}|\mathrm{Re}\,\nu_{0,j}|$.\\
Since $\Psi_z$ depends analytically on $z$, we can apply Stein's
interpolation theorem (\cite[Theorem 4.1]{s-w}), and obtain that,
for every $r'\geq 2+\delta$,
\begin{equation}\label{eq13}
\begin{array}{c}
\Vert \Psi_z \Vert_{L^{q'}(G)}\leq
C_\delta\Vert\xi\Vert_{{p}_0'}\Vert \eta\Vert_{{p}_0'}\textrm{,\quad if }|\mathrm{Re}\, z|\leq \theta_1
\textrm{ and }q':=\frac{r'}{1-|\mathrm{Re}\, z|/\theta_1}.
\end{array}
\end{equation}
But, since $p_0' > 2$, we can choose $\delta >0 $ and $\varepsilon>0$ so small that 
$(1-\frac{\varepsilon}{\theta_1})p_0' \geq 2+\delta$. Then, for $|\mathrm{Re}\, z|\leq \varepsilon$, if we
choose $r':= p_0'(1- \frac{|\mathrm{Re}\, z|}{\theta_1})$ in
\eqref{eq13}, we have $r'\geq 2+\delta$, and hence:
$$\Vert\Psi_z\Vert_{L^{p_0'}(G)}\leq C_\delta\Vert\xi\Vert_{{p}_0'}\Vert\eta\Vert_{{p}_0'}.$$
\end{proof}

\subsection{Proof of Theorem \ref{maintheorem1}}\label{sect6}

Let $p\in [1,\infty[$, $p\neq 2$. The aim is to find a
non-isolated point $\lambda_0$ in the $L^2$-spectrum $\sigma_2(L)$
of $L$ and an open neighbourhood $\mathcal{U}$ of $\lambda_0$ in
$\C$ such that, if $F_0\in C_\infty(\R)$ is an $L^{p}$-multiplier
for $L$, then $F_0$ extends holomorphically to $\mathcal{U}$.
Recall that $C_\infty(\R)$ denotes the space of continuous functions on
$\R$ vanishing at infinity.\\
Since the $L^2$-spectrum of $L$ is contained in $[0,+\infty[$, we
may assume that $F_0\in C_\infty([0,+\infty[)$. Moreover,
according to \cite[Lemma 6.1]{hlm}, it
suffices to consider the case where $2<p'<\infty$.\\
As in the introduction, we can replace $F_0$ by the function
$F=F_0 e^{-\cdot}$, so that $F(L)$ acts on the spaces $L^q(G),\, q\in [p,p']$ by convolution with the function 
$F(L)\delta\in \bigcap_{q=p}^{p'} L^q(G)$. The Kunze-Stein phenomenon implies now that
every $L^p$ function defines a bounded operator on $L^2(G)$ and
also on every Hilbert space $\cal H$ of any unitary representation
$\pi$ of $G$, which is weakly contained in the left regular representation.
Indeed, we know that for any coefficient
$c^\pi_{\xi,\eta}(x):=\langle \pi(x)\xi,\eta\rangle $ of $\pi$, we
have that
$$\Vert c_{\xi,\eta}^\pi\Vert_{p'}\leq C_{p}
\Vert\xi\Vert \Vert \eta\Vert,\ \quad \xi,\eta\in {\cal H},$$
for some constant $C_p>0$. Hence for $f\in L^p(G)$,
$$\vert \int_G f(x) c_{\xi,\eta}^\pi(x) dx\vert\leq \Vert
f\Vert_p \Vert c^\pi_{\xi,\eta}\Vert_{p'}\leq C_p \Vert f\Vert_p
\Vert \xi\Vert \Vert \eta\Vert .$$ Hence there exists a unique
bounded operator $\pi(f)$ on ${\cal H}$, such that $\Vert
\pi(f)\Vert_{\rm op}\leq C_p \Vert f\Vert_p$ and
$$\langle\pi(f)\xi,\eta\rangle=\int_G f(x) c_{\xi,\eta}^\pi(x) dx,
\,\, \ \xi,\eta\in {\cal H}.$$
Choosing now a sequence $(f_\nu)_\nu$ of continuous functions with
compact support, which converges in the $L^p$-norm to
$F(L)\delta$, we see that the operators $\lambda (f_\nu)$ converge
in the operator norm to $\lambda (F(L))=F(\lambda (L)),$ and so
for every unitary representation $(\pi,{\cal H})$ of $G$ which is
weakly contained in the left regular representation $\lambda$, we
have that:
\begin{eqnarray*}
&& \int_G (F(L)\delta)(x) c^\pi_{\xi,\eta}(x)\, dx=
\lim_{\nu\to\infty}\int_{G}f_\nu(x)c^\pi_{\xi,\eta}(x)\, dx\\
&& =\lim_{\nu\to\infty}\langle \pi(f_\nu)\xi,\eta\rangle= \langle
\pi(F(L))\xi,\eta\rangle=\langle F(\pi(L))\xi,\eta\rangle,\qquad
\xi,\eta \in {\cal H}.
\end{eqnarray*}

In particular,
\begin{eqnarray}\label{coefficient}
(F(L)\delta)*(c^\pi_{\xi,\eta})\check{}\,(x) &&=\int_G
(F(L)\delta)(y)\langle \pi(y)\xi,\pi(x)\eta\rangle\,
dy\nonumber\\
&&=\langle F(\pi(L))\xi,\pi(x)\eta\rangle,\qquad \ x\in G,\,
\xi,\eta \in {\cal H}.
\end{eqnarray}

In a first step in order to find $\lambda_0\in\R$ and its
neighborhood ${\cal U}$, we choose a suitable direction $\nu_0$ in
$\a^*$. To this end, let $\omega$ be the Casimir operator of $G$,
and let $\nu\in i\a^*$. Then $\pi_\nu$ is a unitary
representation, and we can define the operator $d\pi_\nu(\omega)$
on the space of smooth vectors in $L^{2}(K/M)$ with respect to
$\pi_\nu$. Moreover, $\pi_{\nu}$ is irreducible (see
\cite[Theorem 1]{kostant}), and therefore
$$d\pi_\nu(\omega) = \chi(\nu)\mathrm{Id},$$
where $\chi$ is a polynomial function on $\a^*$, given by the
Harish-Chandra
isomorphism. Thus, $p$ is in fact a quadratic form.\\
Choose $\nu_0\in \a^*$, $\nu_0\neq 0$, such that $p(\nu_0) \neq 0$.
Then, clearly,
\begin{equation}\label{eq14}
|\chi(iy\nu_0)|\to +\infty\textrm{ as }y\to +\infty\textrm{ in }\R.
\end{equation}
Put $p_0 :=p'$. According to Proposition \ref{interpol}, there is
an $\varepsilon>0$ and a constant $C>0$, such that \eqref{eq9}
holds, for every $z\in U_1:=\{z\in\C|\,|\mathrm{Re}\,z|<\varepsilon\}$. Put
$$\pi_{(z)} := \pi_{z\nu_0}\textrm{ and }\widetilde{T}(z):= T(z\nu_0).$$
Then $(\widetilde{T}(z))_{z\in U_1}$ is an analytic family of
compact
operators on $L^{{p}_0}(K/M)$ (see Proposition \ref{compactop}).\\
And, by an obvious analogue to \cite[Proposition 5.4]{hlm}, there
exists an open connected neighbourhood $U_{y_0}$ of some point
$iy_0$ in $U_1$, with $y_0\in\R$, and two holomorphic mappings
$$\lambda : U_{y_0}\to \C\textrm{ and }\xi: U_{y_0}\to L^{{p}_0'}(K/M)$$
such that, for all $z\in U_{y_0}$ and some constant $C>0,$
\begin{equation}\label{eq6.1}
\begin{array}{ll}
& \widetilde{T}(z)\xi(z) = \lambda(z)\xi(z);\\
& \xi(z)\neq 0 \textrm{ and } \Vert \xi(z)\Vert_{{p}_0'}\leq C.
\end{array}
\end{equation}
Since $\pi_{(iy)}$ is unitary for every $y\in \R$, $\lambda$ is
real-valued on $U_{y_0}\cap i\R$.\\
Fix a non-trivial function $\eta$ in $C^\infty(K/M)$.\\
Let $\Phi_z$, $z\in U_{y_0}$, denote the function on $G$ given by
$$\Phi_z(g) := \langle \pi_{(z)}(g^{-1})\xi(z),\eta\rangle.$$ Then
$\Phi_z(g)$ depends continuously on $z$ and $g$. Moreover, by
\eqref{eq9} and \eqref{eq6.1}, there exists a constant $C_0 >0,$
such that:
\begin{equation}\label{uniform}
\Vert \Phi_z\Vert_{L^{p_0'}(G)}\leq C_0 \textrm{,\quad for all }z\in U_{y_0}.
\end{equation}
Thus, for any $z\in U_{y_0}$, $\Phi_z\in L^{p_0'}(G)$, and
consequently $F(L)\Phi_z\in L^{p_0'}(G)$ is well-defined, since
$F$ is an
$L^{p_0'}$-multiplier for $L$.\\
Put $\mu(z):= -\log \lambda(z)$ ($z\in U_{y_0}$), where $\log$
denotes the principal branch of the logarithm. Since, for $z\in U_{y_0}$, $\xi(z)$ is an eigenvector of
$\widetilde{T}(z)=\pi_{(z)}(h_1)$ associated to the eigenvalue
$\lambda(z)$, where $h_1$ is the convolution kernel of $e^{-L}$,
one has by \eqref{coefficient}, for all $z\in U_{y_0}\cap i\R,\,g\in G:$
\begin{equation}\label{eq6.2}
\begin{array}{rcl}
(F(L)\Phi_z)(g) & = & \langle F(\pi_{(z)}(L))\xi(z), \pi_{(z)}\eta
\rangle\\
& = & F(\mu(z)) \langle \pi_{(z)}(g^{-1})\xi(z),\eta\rangle.
\end{array}
\end{equation}
Let $\psi$ be a fixed element of $C_0(G)$ such that:
$$\int_G \Phi_{iy_0}(x)\psi(x)\, dx \neq 0.$$
By shrinking $U_{y_0}$, if necessary, we may assume that
$\int_G \Phi_{z}(x)\psi(x)\, dx \neq 0$ for all $z\in U_{y_0}$.\\
Then, \eqref{eq6.2} implies that:
\begin{equation}\label{eq6.3}
(F\circ \mu)(z) = \frac{\int_G (F(L)\Phi_z)(x)\psi(x)\,dx} {\int_G
\Phi_{z}(x)\psi(x)\, dx}\textrm{,\quad for }z\in U_{y_0}\cap i\R.
\end{equation}
Observe that the enumerator and the denominator in the right-hand
side of \eqref{eq6.3} are holomorphic functions of $z\in U_{y_0}.$
Indeed, $\langle F(L)\Phi_z,\overline \psi\rangle=\langle \Phi_z, F(L)^*\overline \psi\rangle),$ where 
$F(L)^*\overline\psi\in L^{p_0}$ and $||\Phi_z||_{L^{p_0'}}\le C,$ by \eqref{uniform}.
This implies that the mapping $z\mapsto \langle F(L)\Phi_z,\overline \psi\rangle$ is continuous, and the
holomorphy of this mapping then follows easily from Fubini's and
Morera's
theorems. \\
Therefore, $F\circ\mu$ has a holomorphic extension to
$U_{y_0}.$ \\
Moreover, since $\omega h_1\in L^1(G)$, in view of Proposition
\ref{proprep}, the norm $$\Vert \pi_{(iy)}(\omega
h_1)\Vert_{op}\leq \Vert \omega h_1\Vert_{L^1(G)}$$ is uniformly
bounded, for $y\in \R$. On the other hand, $\pi_{(iy)}(\omega h_1)= d\pi_{(iy)}(\omega)\pi_{(iy)}(h_1) = \chi(iy\nu_0)\pi_{(iy)}(h_1)$,
and so \eqref{eq14} implies that
\begin{displaymath}
\lim_{y\to +\infty}\Vert \widetilde{T}(iy)\Vert =\lim_{y\to
+\infty}\Vert \pi_{(iy)}(h_1)\Vert = 0.
\end{displaymath}
This shows that $\lambda$ is not constant, and hence, varying
$y_0$ slightly, if necessary, we may assume that $\mu'(iy_0)\neq 0$. Then $\mu$ is a local bi-holomorphism near $iy_0$, which
implies, in combination with \eqref{eq6.3}, that $F$ has a
holomorphic extension to a complex neighbourhood of $\lambda_0 :=\mu(iy_0)\in \R$.

\def\ind#1#2{\hb{ind}_{#1}^{#2}}

\section{Transference for $p$-induced representations}
\setcounter{equation}{0}
\subsection{$p$-induced representations}\label{induced}

Let $G$ be a separable locally compact group and $S<G$ a closed
subgroup.\par By \cite{mackey}, there exists a Borel measurable
cross-section $\sigma :G/S\rightarrow G$ for the homogeneous space $H:=G/S$ 
(i.e.~$\sigma (x)\in x$ for every $x\in G/S$) such that $\sigma (K)$ 
is relatively compact for any compact subset $K$ of $H$. Then, every $g\in G$ can be uniquely decomposed as
\begin{equation*}
g=\sigma (x)s,\quad \mbox{with}\; x\in H, s\in S.
\end{equation*}

We put $\Phi :H\times S\rightarrow G,\  \Phi (x,s):=\sigma(x)s$.\par 
Then $\Phi$ is a Borel isomorphism, and we write
\begin{equation*}
\Phi^{-1} (g)=: (\eta (g),\tau (g)).
\end{equation*}

Then
\begin{equation*}
g=\sigma\circ\eta (g)\tau (g),\quad g\in G.
\end{equation*}

For later use, we also define
\begin{eqnarray*}
\tau(g,x) & := & \tau(g^{-1}\sigma(x)),\\
\eta (g,x) & := & \eta (g^{-1}\sigma (x)), \quad g\in G, x\in H\,.
\end{eqnarray*}

Let $dg$ denote the left-invariant Haar measure on $G$,  and
$\Delta_G$ the modular function on $G$, i.e.
\begin{equation*}
\int_G f(gh) dg = \Delta_G (h)^{-1}\int_G f(g)dg,\ h\in G\, .
\end{equation*}

Similarly, $ds$ denotes the left-invariant Haar measure on $S$,
and $\Delta_S$ its modular function.

On a locally compact measure space $Z,$ we denote by $\M_b(Z)$
  the space of all essentially
  bounded measurable functions from $M$ to $\C,$ and by  $\M_0(Z)$
  the subspace of all functions
which have compact support, in the sense that they vanish 
a.e.~outside a compact subset of $Z.$ For $f\in\M_0(G)$, let $\tilde f$
be the function on $G$ given by
\begin{equation*}
\tilde f(g):= \int_S f(gs)\Delta_{G,S} (s)ds,\quad g\in G\, ,
\end{equation*}
where we have put $\Delta_{G,S}(s)=\Delta_G(s)/\Delta_S(s)$, $s\in S$. Then
  $\tilde f$ lies in the space
\begin{eqnarray*}
\E (G,S):=\{\tilde h\in \M_b(G):\tilde h\  \mbox{has compact
support modulo }
S, \mbox{ and }\\
\tilde h (gs)=(\Delta_{G,S}(s))\inv\tilde h(g)\ \mbox{for all}\
g\in G, s\in S  \}
\end{eqnarray*}

In fact, one can show that $\E (G,S)=\{{\tilde f}:f\in \M_0(G)\}$.
Moreover, one checks easily, by means  of the use of a Bruhat function, that 
${\tilde f}=0$ implies $\int_G f(g)dg=0$.\par 
From here it follows that there exists a unique positive linear
functional, denoted by $\int_{G/S}d{\dot g}$, on the space $\E
(G,S)$, which is left-invariant under $G$, such that
\begin{eqnarray}\label{intGS}
\int_G f(g)dg & = & \int_{G/S}\tilde f(g)d{\dot g}\\ \nonumber
 & = & \int_{G/S}\int_S f(gs)\Delta_{G,S}(s)ds\; d{\dot g}\, .
\end{eqnarray}
By means of the cross-section ${\si}$, we can next identify the
function $\tilde h\in\E (G,S)$ with the measurable function $h\in\M_0(H),$ given by
\begin{equation*}
h(x):=R\tilde h(x):=\tilde h (\si (x)),\quad x\in H\, .
\end{equation*}
Notice that, given $h\in\M_0(H),$ the corresponding function
$\tilde h=:R\inv h\in\E (G,S)$ is given by
\begin{equation*}
\tilde h(\si (x)s)=h(x)\Delta_{G,S}(s)\inv\, .
\end{equation*}

The mapping $h\mapsto\int_{G/S}\tilde h (g)d{\dot g}$ is then a
positive Radon measure on $C_0(H)$, so that there exists a unique
regular Borel measure $dx$ on $H=G/S$, such that
\begin{equation}\label{Bomeas}
\int_{G/S}\tilde h(g)d{\dot g}=\int_H h(x)dx,\quad h\in C_0(H)\, .
\end{equation}

Formula \eqref{intGS} can then be re-written as
\begin{equation}\label{intHintS}
\int_G f(g) dg=\int_H\int_S f(\si(x)s)\Delta_{G,S}(s)ds\; dx\, .
\end{equation}
Notice that the left-invariance of $\int_{G/S}d{\dot g}$ then
translates into the following quasi-invariance property of the
measure $dx$ on $H$:
\begin{equation}\label{dxH}
\int_H h(\eta (g,x))\Delta_{G,S}(\tau (g,x))\inv dx=\int_H h(x)dx\ \mbox{ for every}\; g\in G\, .
\end{equation}
Formula \eqref{intHintS} remains valid for all $f\in L^1(G)$.\\
Next, let $\rho$ be a strongly continuous isometric representation
of $S$  on a
  complex Banach space $(X,||\cdot||_X),$ so that in particular
\begin{equation*}
||\rho(s)v||_X=||v||_X\ \mbox{ for every }s\in S, v \in X.
\end{equation*}
Fix $1\le p<\iy,$ and let $L^p(G,X;\, \rho)$ denote the Banach
space of all Borel measurable functions $\stacksx : G\rightarrow X$, which satisfy the covariance condition
\begin{equation*}
\stacksx (gs)=\Delta_{G,S}(s)^{-1/p}\rho (s\inv)[\stacksx (g)],
\quad\mbox{for all}\; g\in G,s\in S\, ,
\end{equation*}
and have finite $L^p$-norm $||\stacksx ||_p :=
\l(\int_{G/S}||\stacksx (g)||^p_{X}d{\dot g}\r)^{1/p}.$\\
Notice that the function $g\mapsto ||\stacksx (g)||^p_{X}$ 
satisfies the covariance property of
functions in $\E(G,S)$, so that the integral
$\int_{G/S}||\stacksx (g)||^p_{X}d{\dot g}$ is well-defined.\\

The {\it $p$-induced representation} $\pi_p={\rm ind}^G_{p,S}\rho$
is then the left-regular representation $\lambda_G=\lambda$ of $G$ acting on $L^p(G,X;\rho)$, i.e.
\begin{equation*}
\left[ \pi_p(g)\stacksx \right](g'):=\stacksx(g\inv g'),\quad g,g'\in G,\, \tilde \xi\in L^p (G,X;\,\rho).
\end{equation*}
By means of the cross-section $\si$, one can realize $\pi_p$ on
the $L^p$-space $L^p(H,X).$\\
To this end, given $\stacksx\in L^p(G,X;\rho)$, we define $\xi\in L^p(H,X)$ by
\begin{equation*}
\xi (x):=\T\stacksx (x):=\stacksx (\si (x)),\quad  x\in H\, .
\end{equation*}
Because of \eqref{Bomeas}, $\T:L^p(G,X;\,\rho)\rightarrow L^p(H,X)$ is a linear isometry, with inverse
\begin{equation*}
\T\inv\xi (\si(x)s):=\stacksx (\si(x)s)=\Delta_{G,S}(s)^{-1/p}\rho
(s\inv) [\xi (x)]\, .
\end{equation*}
Since, for $g\in G,\, y \in H$ and $\stacksx\in L^p(G,X;\, \rho)$,
\begin{eqnarray*}
\stacksx (g\inv\si (y)) & = & \stacksx (\si\circ\eta (g\inv\si
(y))\tau (g\inv\si (y)))\\
 & = &  \stacksx (\si (\eta (g,y))\tau(g,y))\\
 & = & \Delta_{G,S}(\tau (g,y))^{-1/p} \rho (\tau (g,y)\inv)
\left[\stacksx (\si (\eta(g,y)))\right],
\end{eqnarray*}
we see that the induced representation $\pi_p$ can also be
realized on $L^p(H,X)$, by
\begin{eqnarray}\label{pip}
\left[ \pi_p(g)\xi\right] (y) = \Delta_{G,S} (\tau
(g,y))^{-1/p}\rho (\tau (g,y)\inv)\left[ \xi (\eta (g,y))\right],
  \end{eqnarray}
for  $ g\in G,y\in H,\; \xi\in L^p (H,X)\, .$

Observe that $\pi_p(g)$ acts isometrically on $L^p (H,X)$, for
every $g\in G$. This is immediate from the original realization of
$\pi_p$ on $L^p(G, X;\, \rho)$, but follows also from \eqref{dxH}, 
in the second realization given by \eqref{pip}.
\begin{exs}\label{example1}${}$
{\rm
\begin{itemize}
\item[{\bf (a)}] If $S\lhd G$ is a closed, normal subgroup, then $H=G/S$ is
again  a group, and one finds that, for a suitable normalization
of the left-invariant Haar measure $dx$ on $H$, we have
\begin{equation*}
\int_G f(g)dg = \int_H\int_S f(\si (x)s)ds\, dx,\quad  f\in L^1(G)\, .
\end{equation*}
In particular, $\Delta_G | _S=\Delta_S$, so that $\Delta_{G,S}=1$
and $dx$ in
\eqref{intHintS} agrees with the left-invariant Haar measure on $H$.\\
Furthermore, there exists a measurable mapping $q:H\times H\rightarrow S$\, ,\par such that
\begin{equation*}
\si (x)\inv\si (y) = \si (x\inv y)q(x,y),\quad x,y\in H\, ,
\end{equation*}
since $\si (x)\inv\si (y)\equiv \sigma(x\inv y)$ modulo $S$. Thus, if
$g=\si (x)s$, then
\begin{eqnarray*}
g\inv\si (y) & = & s\inv\si (x)\inv\si (y) = s\inv\si (x\inv y)q(x,y)\\
& = & \si (x\inv y)((s\inv)^{\si (x\inv
y)\inv}q(x,y)).
\end{eqnarray*}
(Here we use the notation $ s^g:= g s g\inv ,\ s\in S, g\in G$.)\\
This shows that $\tau (g,y) = (s\inv)^{\si (x\inv  y)\inv}q(x,y)$ and 
$\eta (g,y) = x\inv y$. 
Hence $\pi_p$ is given as follows:
\begin{eqnarray}\label{pips}
  \left[ \pi_p (\si (x)s)\xi\right](y)
  =  \rho (q (x,y)\inv s^{\si (x\inv  y)\inv})\left[ \xi (x\inv y)\right],
\end{eqnarray}
for $(x,s)\in H\times S,\; y\in H,\ \xi\in L^p (H,X).$\\
We remark that it is easy to check that:
\begin{equation*}
q(x,y)\inv s^{\si (x\inv  y)\inv} = s^{\si (y)\inv\si(x)}q(x,y)\inv .
\end{equation*}
Notice that \eqref{pips} does not depend on $p$.
\item[{\bf (b)}] In the special case where $\rho =1$ and $S$ is normal, the
induced
  representation  $\iota={\rm ind}_S^G 1$ is given by
\begin{equation*}
\left[ \iota(\si (x)s)\xi\right](y)=\xi (x\inv y)\, .
\end{equation*}
For the integrated representation, we then have
\begin{eqnarray*}
\left[ \iota(f)\xi\right](y) & = & \int_H\int_S f(\si(x)s)\xi(x\inv y)ds\, dx\\
 & = & \int_H \tilde f (x)\xi (x\inv y) dx\\
 & = & \left[ \lambda_H (\tilde f)\xi\right](y)\, ,
\end{eqnarray*}
i.e.
\begin{equation*}
\iota(f) = \lambda_H (\tilde f)\, ,
\end{equation*}
where
\begin{equation*}
\tilde f (x) := \int_S f(\si (x)s)ds \, ,
\end{equation*}
i.e.~$\tilde f$ is the image of $f$ under the quotient map from
$G$ onto $G/S$.
\end{itemize}
}
\end{exs}


\subsection{A transference principle}\label{transfer}
If $\xi \in L^p(H,X)$, and if $\phi :S\rightarrow\C$, we define
the $\rho$-{\it twisted tensor product}
\begin{displaymath}
\begin{array}{c}
\xi\otimes_{\rho}^p\phi : G\rightarrow X\quad\mbox{by}\\
\left[\xi\otimes_{\rho}^p\phi\right] (\si (x)s):=\phi(s)\Delta_{G,S}(s)^{-1/p}\rho (s\inv)[\xi (x)],\quad (x,s)\in H\times S\, .
\end{array}
\end{displaymath}
Let us denote by $X^*$ the dual space of $X.$ For any complex
vector space $Y,$ we denote by $\overline Y$ its {\rm complex conjugate,} which, as an additive group, is the space $Y,$ but
with  scalar multiplication given by $\overline \lambda y,$ for
$\lambda\in \C$ and $y\in Y.$ In the sequel, we assume
  that
$X$ contains a dense, $\rho$- invariant subspace $X_0,$ which
embeds via an anti-linear mapping $i:X_0\hookrightarrow\overline{X^*}$ into the complex conjugate of the dual space of
$X,$ in such a way that
\begin{equation}\label{normid}
||x||=\sup_{\{v\in X_0: \,||i(v)||_{X^*}=1\}}|\lan x,v\ran |\quad
\mbox{for every } x\in X.
\end{equation}
Here, we have put
$$\lan x,v\ran :=i(v)(x),\quad v\in X_0,x\in X.$$
Moreover, we assume that
\begin{equation}\label{normid*}
||i(\rho(s)v)||_{X^*}=||i(v)||_{X^*}\quad \mbox{for every }v\in X_0, s\in S,
\end{equation}
and
\begin{equation} \label{iso*}
\lan \rho(s)x,\rho(s)v\ran=\lan x,v\ran \quad \mbox{for every }x\in X,v\in X_0,s\in S.
\end{equation}

The most important example for us will be an $L^p$-space
$X=L^p(\Om),\ 1\le p<\infty,$ on a measure space
$(\Omega,d\omega),$ and a representation $\rho$ of $G$ which acts
isometrically on $L^p(\Om)$ as well as on its dual space
$L^{p'}(\Om)$ (i.e.~$\frac{1}{p}+\frac{1}{p'}=1$). In this case,
by interpolation, we have $||\rho(g)\xi||_r\le ||\xi||_r$, for
$|\frac 1r-\frac 12|\le|\frac 1p-\frac 12|,\ g\in G,$ which
implies that indeed $\rho(g)$ acts isometrically on $L^r(\Om),$
for $|\frac 1r-\frac 12|\le|\frac 1p-\frac 12|.$ In particular,
$\rho$ is a unitary representation on 
$L^2(\Om).$ We can then choose $X_0:=L^{p'}(\Om)\cap L^{p}(\Om)\subset L^2(\Om),$ and put
$$
i(\eta)(\xi):=\int_\Om \xi(\omega)\overline{\eta(\omega)}\, d\omega,\quad
\eta\in L^{p'}(\Omega)\cap L^{p}(\Omega), \xi\in L^{p}(\Omega).
$$
Notice that \eqref{normid*} and \eqref{iso*} are always satisfied,
if $\rho$ is a unitary character.
\begin{lemme}\label{lem}
Let $\phi\in L^p(S),\psi\in L^{p'}(S),\xi\in L^p(H,X_0)$ and
$\eta\in L^{p'}(H,X_0)$, where $\frac{1}{p}+\frac{1}{p'}=1$. Then,
for every $g\in\, G,$
\begin{eqnarray}\label{identt}
\lan\lambda_G (g)\left(\xi\otimes_{\rho}^p\phi\right) ,\,
\eta\otimes^{p'}_{\rho} \psi\ran = \int_H \phi\ast
\stackrel{\veebar}{\psi}  (\tau (g,x))\lan [\pi_p(g)\xi] (x),\,
\eta (x)\ran\, dx.
\end{eqnarray}
\end{lemme}
\begin{proof} By \eqref{intHintS}, we have
\begin{eqnarray*}
&&\lan\lambda_G (g)\left(\xiphi\right) ,\, \etapsi\ran\\
  =&&  \int_H\int_S\lan \xiphi (g^{-1}\si (x)s),\, \etapsi
(\si(x)s)\ran\, \Delta_{G,S} (s) ds\, dx\\
  = &&\int_H\int_S\lan \xiphi \left(\si (\eta (g,x))\tau (g,x)s\right),\,
\etapsi (\si (x)s)\ran\, \Delta_{G,S} (s) ds\, dx\\
  = &&\int_H\int_S\Delta_{G,S} (\tau (g,x)s)^{-\frac{1}{p}}
\Delta_{G,S} (s)^{-\frac{1}{p'}}\phi (\tau (g,x)s)\ol{\psi}(s)\\
&& \hskip2cm\lan \rho (s\inv\tau (g,x)\inv)\, [\xi (\eta
(g,x))],\, \rho (s\inv)[\eta (x)]\ran
\,\Delta_{G,S}(s)\, ds\, dx\\
  = && \int_H\int_S\Delta_{G,S}
  (\tau (g,x))^{-\frac{1}{p}}\phi (\tau (g,x)s)\ol{\psi}(s)\, ds\\
&& \hskip2cm\lan \rho (\tau (g,x)\inv)[\xi (\eta (g,x))],\, \eta
(x)\ran\, dx\, .
\end{eqnarray*}
Here, we have used that, by \eqref{iso*}, $\lan \rho (s\inv)v_1,\rho (s\inv) v_2\ran=\lan v_1,v_2\ran$ 
for all $v_1,v_2\in\, X_0.$\\
But,
$$
\int_S\phi (\tau (g,x)s)\overline\psi (s) ds
   =   \int_S\phi (s)\psi (\tau (g,x)\inv s) ds
   =   \phi\ast\stackrel{\veebar}{\psi}(\tau (g,x)),
$$
and
$$\Delta_{G,H} (\tau (g,x))^{-\frac{1}{p}}\lan \rho (\tau (g,x)\inv)[\xi
(\eta (g,x))], \eta (x)\ran
  = \lan [\pi_p (g)\xi](x),\eta (x)\ran\, ,
$$
and thus \eqref{identt} follows.
\end{proof}

From now on, we shall assume that the group $S$ is {\sl amenable}.\\

Since $G$ is separable, we can then choose an increasing sequence
$\{A_j \}_j$ of compacta in $S$ such that $A_j\inv =A_j$ and $S=\bigcup\limits_{j}A_j$, and put
\begin{eqnarray*}
\phi_j=\phi_j^p :=\frac{\chi_{A_{j}}}{|\, A_j\, |^{1/p}},\quad
\psi_j=\psi_j^{p'} :=\frac{\chi_{A_{j}}}{|\, A_j\, |^{1/p'}},
\end{eqnarray*}
where $\chi_A$ denotes the characteristic function of the subset
$A$. Then $\stackrel{\veebar}{\psi_j}=\psi_j ,\ || \phi_j ||_p= ||\psi_j ||_{p'}=1$, and, because of the amenability of $S$ (see
\cite{pier}), we have
\begin{eqnarray}\label{to1}
\chi_j :=\phi_j\ast\psi_j\;\mbox{tends to 1, uniformly on compacta in $S$}\, .
\end{eqnarray}

\begin{proposition}\label{ind}
Let $\pi_p={\rm ind}^G_{p,S}\,\rho$ be as before, where $S$ is
amenable, and let $\xi, \eta\in C_0 (H,X_0).$\par
Then
\begin{equation*}
\lan \pi_p (g)\xi ,\eta\ran = \lim_{j\rightarrow\iy}\lan
\lambda_G(g)(\xiphi_j),\etapsi_j\ran\, ,
\end{equation*}
uniformly on compacta in $G$.
\end{proposition}

\begin{proof} By Lemma \ref{lem},
\begin{eqnarray*}
  \lan \lambda_G (g) (\xiphi_j),\etapsi_j\ran
  = \int_H \chi_j (\tau (g,x))\lan [\pi_p(g)\xi](x),\eta (x)\ran \,dx\, .
\end{eqnarray*}
Fix a compact set $K=K\inv\subset H$ containing the supports of
$\xi$ and $\eta$, and let $Q\subset G$ be any compact set. We want to prove that $\{\tau (g,x)|\,g\in\, Q, x\in\, K\}$ is 
relatively compact, for then, by \eqref{to1}, we immediately see that
\begin{eqnarray*}
\lim_{j\rightarrow\iy}\lan \lambda_G(g)(\xiphi_j),\etapsi_j\ran
=\int_H\lan [\pi_p(g)\xi](x),\eta (x)\ran dx =\lan
\pi_p(g)\xi,\eta\ran\, ,
\end{eqnarray*}
uniformly for $g\in\, Q$.\\
Recall that $\tau (g,x)=\tau (g\inv\si (x))$. Therefore, since
$\si (K)$ is relatively compact, it suffices to prove that $\tau$
maps compact subsets of $G$ into relatively compact sets in $S$.
So, let again $Q$ denote a compact subset of $G$, and put
$M:=Q\mbox{ mod }\, S\subset H=G/S$. Then $M$ is compact, so that
$\ol{\si(M)}$ is compact in $S$. And, since $\tau (\si (x)s)=s$
for every $x\in\, H,\, s\in\, S$, we have
\begin{eqnarray*}
\tau (Q) & = & \{s\in\, S|\,\si (x)s\in\, Q\,\mbox{ for some}\,
x\in\, M\}\\
 & = & \si (M)\inv Q\, ,
\end{eqnarray*}
which shows that $\tau (Q)$ is indeed relatively compact.
\end{proof}

\begin{thm}\label{intermediatetheorem}
For every bounded measure $\mu \in\, M^1(G)$, we have
\begin{equation*}
||\,\pi_p (\mu)\,||_{L^p(H,X)\rightarrow L^p(H,X)}\le
||\,\lambda_G(\mu)\,||_{L^p(G,X) \rightarrow L^p(G,X)}\, .
\end{equation*}
\end{thm}
\begin{proof} Let $\xi,\eta\in\, C_0 (H,X_0).$ Observe first that, for $g\in G,$
\begin{eqnarray*}
&&|\, \lan \lambda_G(g)(\xiphi_j),\etapsi_j\ran\, |\\
  \le && ||\, \lambda_G(g)(\xiphi_j)\, ||_{L^p(G,X)}\ ||\, i\circ(\etapsi_j)\, ||_{L^{p'}(G,X^*)}\\
  =   &&  ||\, \xiphi_j\,||_{L^p(G,X)}\  ||\,i\circ(\etapsi_j)\,
||_{L^{p'}(G,X^*)}\, ,
\end{eqnarray*}
where
\begin{eqnarray*}
&& ||\, \xiphi_j\,||_{L^p(G,X)}^p\\
& = & \int_H\int_S |\, \phi_j(s)\, |^p \Delta_{G,S}(s)\inv ||\,\rho (s\inv)[\xi (x)]\, ||^p_{X}\Delta_{G,S}(s) ds\, dx\\
  & = & \int_S |\, \phi_j(s)\, |^p ds\quad \int_H ||\, \xi (x)\,||^p_{X} dx\\
& = & ||\, \xi\, ||^p_{L^p(H,X)}\, ,
\end{eqnarray*}
since $\rho (s\inv)$ is isometric on $X$, so that
\begin{equation}\label{id1}
||\, \xiphi_j\,||_{L^p(G,X)} = ||\xi ||_{L^p(H,X)},
\end{equation}
and similarly, because of \eqref{normid*},
\begin{equation}\label{id2}
||\,i\circ(\etapsi_j)\, ||_{L^{p'}(G,X^*)}=||\, i\circ\eta\,||_{L^{p'}(H, X^*)}\, .
\end{equation}
This implies
\begin{eqnarray*}
|\, \lan \lambda_G(g)(\xiphi_j),\etapsi_j\ran |\;
\le  ||\xi ||_{L^p(H,X)}\; ||\, i\circ\eta\, ||_{L^{p'}(H, X^*)} .
\end{eqnarray*}
Therefore, if $\mu\in\, M^1(G)$, Proposition \ref{ind} implies, by
the dominated convergence theorem, that
\begin{equation}\label{lim}
\lan \pi_p(\mu)\xi, \eta\ran = \lim_{j\rightarrow\iy}\lan\lambda_G(\mu)(\xiphi_j), \etapsi_j\ran\, .
\end{equation}
Moreover, by \eqref{id1} and  \eqref{id2},
\begin{eqnarray*}
&&|\, \lan \lambda_G(\mu)(\xiphi_j),\etapsi_j\ran\, |\nonumber\\
  &&= |\, \int_G\lan\lambda_G(\mu) (\xiphi_j),\, \etapsi_j\ran\, dg |\nonumber\\
&&\le  ||\,\lambda_G(\mu)\,||_{L^p(G,X) \rightarrow L^p(G,X)}\,||\xi ||_{L^p(H,X)}\; ||\, i\circ\eta\, ||_{L^{p'}(H, X^*)} .
\end{eqnarray*}
By \eqref{lim}, we therefore obtain
\begin{equation}\label{xinu}
|\, \lan \pi_p(\mu)\xi ,\eta\ran\, |\le
||\,\lambda_G(\mu)\,||_{L^p(G,X) \rightarrow L^p(G,X)}\,  ||\xi ||_{L^p(H,X)}\; ||\, i\circ\eta\, ||_{L^{p'}(H, X^*)} .
\end{equation}
In view of \eqref{normid}, this implies the theorem, since $C_0(H,X_0)$ lies dense in $L^p(H,X).$
\end{proof}

\begin{cor1}[Transference]\label{trans}
Let $X=L^p(\Omega).$ Then,  for every $\mu \in\, M^1 (G)$, we have
\begin{equation*}
||\,\pi_p (\mu)\,||_{L^p(H,L^p(\Omega))\rightarrow
L^p(H,L^p(\Omega)))}\le ||\,\lambda_G(\mu)\,||_{L^p(G)\rightarrow L^p(G)}\, .
\end{equation*}
\end{cor1}

\begin{proof} If $X=L^p(\Omega)$ and $h\in L^p(G,X),$ then, by Fubini's theorem,
$$
||\lambda_G(\mu)h||_{L^p(G,X)}^p =\int_\Omega||\mu\ast h(\cdot,\omega)||_{L^p(G)}^p\, d\omega 
\le ||\,\lambda_G(\mu)\,||_{L^p(G)\rightarrow L^p(G)}^p||h||_{L^p(G,X)},
$$ 
hence
$$
||\lambda_G(\mu)||_{L^p(G,X)\rightarrow L^p(G,X)}
\le ||\,\lambda_G(\mu)\,||_{L^p(G)\rightarrow L^p(G)}.
$$ 
In combination with \eqref{xinu}, we obtain the desired estimate.
\end{proof}

\begin{rem}
{\rm 
We call a Banach space $X$ to be of {\bf $L^p$-type},
$1\le p< \infty,$ if there exists an embedding $\iota:X\hookrightarrow L^p(\Omega)$ into an
$L^p$-space such that
$$
\frac 1C ||x||_X\le ||\iota(x)||_{L^p(\Omega)}\le C||x||_X\quad\mbox{ for every } x\in X,
$$ 
for some constant $C\ge 1.$\\
For instance, any separable Hilbert space $\mathcal H$ is of
$L^p$-type, for $1\le p < \infty,$ or, more generally, any space
$L^p(Y,\mathcal H).$ This follows easily from Khintchin's
inequality. Corollary \ref{trans} remains valid for spaces $X$ of
$L^p$-type, by an obvious modification of the proof.
}
\end{rem}
Denote by $C^{\ast}_r(G)$ the reduced $C^{\ast}$-algebra of $G$.
If $p=2$, we can extend \eqref{lim} to $C^{\ast}_r(G)$.
\begin{proposition}\label{limite}
If $p=2$ and $X=L^2(\Omega),$ then the unitary representation
$\pi_2$ is weakly contained in the left-regular
  representation $\lambda_G$. In particular, for any $K\in C^{\ast}_r(G)$,
the operator
  $\pi_2(K)\in\, {\cal B}(L^2(H,L^2(\Omega)))$ is well-defined.\par
Moreover, for all $\xi ,\eta\in\, C_0(H,L^2(\Omega))$, we have
\begin{equation}\label{pi2K}
\lan \pi_2(K)\xi ,\eta\ran\quad =\quad\lim_{j\rightarrow\iy}\lan
\lambda_G(K)
(\xi\otimes^2_{\rho}\phi_j),\eta\otimes^2_{\rho}\phi_j,\ran\, .
\end{equation}
\end{proposition}

\begin{proof} If $K\in\, C^{\ast}_r(G)$, then we can find a sequence
$\{f_k\}_k$ in $L^1(G)$, such that
$\lambda_G(K)=\displaystyle{\lim_{k\rightarrow\iy}}\lambda_G(f_k)$
in the operator norm $||\, \cdot\, ||$ on $L^2(G)$. But,
\eqref{xinu} implies that
\begin{equation}\label{pi2}
||\, \pi_2(f)\, ||\le ||\, \lambda_G(f)\, ||\, ,\;\mbox{ for
all}\; f\in\, L^1(G),
\end{equation}
where $||\, \cdot\, ||$ denotes the  operator norm on ${\cal B}(L^2(H,L^2(\Omega)))$ and 
${\cal B}(L^2(G)),$ respectively. Therefore, the $\{\pi_2(f_k)\}_k$
form a Cauchy sequence in ${\cal B}(L^2(H,L^2(\Omega))),$ whose limit we denote by $\pi_2(K)$.\\
It does not depend on the approximating sequence $\{f_k\}_k$.
Moreover, from \eqref{pi2} we then deduce that
\begin{equation}\label{laGK}
||\,\pi_2(K)\, ||\,\le\;||\, \lambda_G(K)\, ||\; =\; ||\, K\,
||_{C^{\ast}_r(G)}\, ,\; \mbox{ for all}\; K\in\, C^{\ast}_r(G)\,
.
\end{equation}
In particular, we see that $\pi_2$ is weakly contained in
$\lambda_G$. It remains to show
\eqref{pi2K}.\\
Given $\varepsilon >0$, we choose $f\in\, C_0(G)$ such that $||\, K-f\, ||_{C^{\ast}_r(G)}<\varepsilon/4$. 
Next, by \eqref{xinu}, we can find $j_0$ such that
\begin{equation*}
|\lan \pi_2(f)\xi ,\eta\ran\;
-\;\lan\lambda_G(f)(\xi\otimes^2_{\rho}\phi_j),\,
  \eta\otimes^2_{\rho}\phi_j\ran\,
| < \varepsilon/4\quad\mbox{for all }\, j\ge j_0\, .
\end{equation*}
Assume without loss of generality that $||\, \xi\, ||_2 = ||\,\eta\, ||_2=1$. Then, by \eqref{laGK},
\begin{eqnarray*}
|\, \lan \pi_2(K)\xi, \eta\ran\; -\;\lan \pi_2(f)\xi,\eta\ran\, |
\le   ||\, K-f\, ||_{C^{\ast}_r(G)}\; ||\, \xi\, ||_2\; ||\, \eta
  ||_2\quad <\varepsilon/4,
\end{eqnarray*}
and furthermore
\begin{eqnarray*}
 && |\lan \lambda_G(K)(\xi\otimes^2_{\rho}\phi_j),\eta\otimes^2_{\rho}\phi_j\ran-\lan\lambda_G(f)(\xi\otimes^2_{\rho}\phi_j),
\eta\otimes^2_{\rho}\phi_j\ran\, |\\
&\le & ||\, K-f\, ||_{C^{\ast}_r(G)}\, ||\;
\xi\otimes^2_{\rho}\phi_j\, ||_2\; ||\,
\eta\otimes^2_{\rho}\phi_j\, ||_2\\
& < & \frac{\varepsilon}{4}\; ||\,\xi\, ||_2\; ||\, \eta\, ||\,_2 =\varepsilon/4\, .
\end{eqnarray*}
Combining these estimates, we find that
\begin{eqnarray*}
|\, \lan \pi_2(K)\xi ,\eta\ran -\lan\lambda_G(K)(\xi\otimes^2_{\rho}\phi_j),
\eta\otimes^2_{\rho}\phi_j\ran\, | <\frac{\varepsilon}{4}+\frac{\varepsilon}{4}+\frac{\varepsilon}{4}<\varepsilon
\quad \mbox{for all }\, j\ge j_0\, .
\end{eqnarray*}
\end{proof}

\begin{cor2}\label{passage2} 
Assume that $\rho$ is a unitary representation on a separable Hilbert space $X$, 
for instance a unitary character of $S,$ and that $\Delta_{G,S}=1.$ 
Let $K\in C_r^*(G),$ and assume that $\lambda_G(K)$ extends from $L^2(G)\cap L^p(G)$ to a
bounded linear operator on $L^p(G),$ where $1\le p<\infty.$\\
Then $\pi_2(K)$ extends from $L^2(G/S)\cap L^p(G/S)$ to a bounded
linear operator on $L^p(G/S),$ and
\begin{equation}\label{ne1}
||\pi_2(K)||_{L^p(G/S)\to L^p(G/S)}\le ||\lambda_G(K)||_{L^p(G)\to L^p(G)}.
\end{equation}
Moreover, for $f\in L^1(G),$ we have $\pi_p(f)=\pi_2(f)$ on $C_0(G/S).$
\end{cor2}
\begin{proof}
 If $\xi,\eta\in C_o(H),$ then, since $\Delta_{G,S}=1,$
$$
\lan \pi_2(K)\xi,\eta\ran=\lim_{j\to\infty}\lan \lambda_G(K)(\xi\otimes^2_{\rho}\phi_j^2),
\eta\otimes^2_{\rho}\psi_j^2\ran=\lim_{j\to\infty} \lan
\lambda_G(K)(\xi\otimes_{\rho}^p\phi_j^p),
\eta\otimes^{p'}_{\rho}\psi_j^{p'}\ran.
$$
And,
\begin{eqnarray*}
&&|\lan \lambda_G(K)(\xi\otimes_{\rho}^p\phi_j^p),
\eta\otimes^{p'}_{\rho}\psi_j^{p'}\ran| \nonumber \\
&\le & ||\lambda_G(K)||_{L^p(G)\to L^p(G)} \,
||\xi\otimes_{\rho}^p \phi_j^p||_{L^p(G)}\,
||\eta\otimes^{p'}_{\rho}\psi_j^{p'}||_{L^{p'}(G)}\nonumber\\
&\le & ||\lambda_G(K)||_{L^p(G)\to L^p(G)} \, ||\xi||_{L^p(H)}\,
||\eta||_{L^{p'}(H)}.
\end{eqnarray*}
Estimate \eqref{ne1} follows.\\
That $\pi_p(f)=\pi_2(f)$ on $C_0(H),$ if $f\in L^1(G),$ is
evident, since $\Delta_{G,S}=1.$
\end{proof}
\section{The case of a non-compact semi-simple factor}\label{sectbig4}
\setcounter{equation}{0}
%
In this section, we shall give our proof of Theorem \ref{maintheorem2}.\\ 
Let us first notice the following consequence of Corollary 
\ref{passage2}.\\

Assume that
$S$ is a closed, normal and amenable subgroup of $G,$ and  let $L=-\sum_j X_j^2$  be
a sub-Laplacian on $G.$ Denote by $\iota_2:={\rm ind}_S^G 1$ the representation of $G$ induced by the trivial
character of $S$ (compare Example \ref{example1}), and  let $\tilde{L}=-\sum_j(X_j\text { mod } \s)^2=d\iota_2(L) $ 
be  the corresponding sub-Laplacian on the quotient group $H:=G/S.$ Then
\begin{equation}\label{multtrans}
{\cal M}_p(L)\cap C_\infty (\R)\subset {\cal M}_p(\tilde L)\cap
C_\infty (\R).
\end{equation}
In particular, if $\tilde L$ is of holomorphic $L^p$-type, then so is $ L.$\\
In order to prove \eqref{multtrans}, assume that  $F$ is an
$L^p$-multiplier for $L$
contained in $C_{\infty}(\R).$ Then
$F(L)$ lies in $C^*_r(G),$ and by Corollary \ref{passage2} the operator
$\iota_2(F(L))=F(d\iota_2(L))=F(\tilde L)$ extends from $L^2(H)\cap L^p(H)$ to
a bounded operator on $L^p(H),$ so that $F\in {\cal M}_p(\tilde L)\cap C_\infty (\R).$\\

Let now $G$ be a connected Lie group, with radical $S=\exp\s.$ Then there
exists
a connected, simply connected semi-simple Lie group $H$ such that $G$ is the
semi-direct product of  $H$ and $S,$ and this Levi factor $H$ has a
discrete center
$Z$ (see \cite{bourbaki}). Let $L$ be a sub-Laplacian on $G,$ and denote by
$\tilde L$ the corresponding sub-Laplacian on $G/S\simeq H$ and by
$\tilde{\tilde
L}$ the sub-Laplacian on $H/Z$ corresponding to $\tilde L$ on $H.$ We have 
that $Z$ and $S$
are amenable groups, and $H/Z$ has finite center. From Theorem
\ref{maintheorem1}, we thus find that  $\tilde{\tilde L}$ is of holomorphic
$L^p$-type for every
$p\ne 2,$  if we assume that $H$ is non-compact, and
\eqref{multtrans} then allows us to conclude that the same is true of
$\tilde L,$
and then also of $L.$

\section{Compact extensions of exponential solvable Lie
groups}\label{sectbig5}
\setcounter{equation}{0}
\subsection{Compact operators arizing in induced representations}
\label{compactext}

Let now $K=\exp\k$ be a connected compact Lie group acting
continuously on an exponential solvable Lie group $S=\exp{\s}$ by
automorphisms $\sigma(k)\in \hb{ Aut (S)}, k\in K$. We form the
semi-direct product $G=K\ltimes S$ with the multiplication given
by:
$$(k,s)\cdot (k',s')=(k k',\sigma({k'}^{-1}) s s'),\ \ k,k'\in K, s,s'\in
S.$$ The left Haar measure $dg$ is the product of the Haar measure
of $K$ and the left Haar measure of $S$. Let us choose a
$K$-invariant scalar product $\sp{\cdot}{\cdot}$ on the Lie algebra $\s $ of
$S$. Denote by $\n$ the nil-radical of $\s$. Since every
derivation $d$ of $\s$ maps the vector space $\s$ into the
nil-radical, it follows that the orthogonal complement $\b$ of
$\n$ in $\s$ is in the kernel of $d\sigma(X)$ for every $X\in \k$.
The following decomposition of the solvable Lie algebra $\s$ has
been given in \cite{boidol}. Choose an element $X\in\b$, which is
in general position for the roots of $\s$, i.e., for which
$\lambda(X)\ne \mu(X)$ for all roots $\mu\ne\lambda$ of $\s$. Let
$\s_0=\pa{ Y\in\s;\ \hb{ad}^l(X)Y=0 \hb{ for some }l\in\N^*}$.
Then $\s_0$ is a nilpotent subalgebra of $\s$, which is
$K$-invariant  (since $[X,\k]=\pa 0$) and $\s=\s_0+\n$. Let $\a$
be the orthogonal complement of $\n\cap \s_0$ in $\s_0$. Then $\a$
is also a $K$-invariant subspace of $\s$ (but not in general a
subalgebra) and $\s=\a\oplus \n$. Let $N=\exp{\n}\subset S$ be the
nil-radical of the group $S$. Then $S$ is the topological product
of $A=\exp\a$ and $N$. Finally our group $G$ is the topological
product of $K,A$ and $N$. Hence every element $g$ of $G$ has the
unique decomposition:
$$g=k_g\cd a_g\cd n_g, \hb{ where } k_g\in K,\ a_g\in A \hb{ and }
n_g\in N.$$
We shall use the notations and constructions of \cite{hlm} in the
following but we have to replace there the symbol $G$ with the
letter $S$.\\
%
Let $h:G\to\C$ be a function. For every $x\in G$, we denote by
$\tilde{h}(x)$ the function on $S$ defined by:
$$\tilde{h}(x)(s)=h(xs),\ s\in S.$$ 
Also for a function $r:S\to \C$ and for $x\in G$, we let ${}^xr:S\to\C$ 
be defined by:
$${}^xr(s):=r(x sx\inv).$$
\\
We say that a Borel measurable function $\o: G\to\R_+ $ is a
weight, if $1\leq\o(x)=\o(x\inv)$ and $\o(x y)\leq \o(x)\o(y)$ for
every $x,y\in G$. Then the space 
$$L^p(G,\o)=\pa{ f\in L^p(G)| \,\no{f}_{\o,p}:= 
\int_G \va{f(g)}^p \o(g)\,dg<\iy},$$
for $1\leq p\leq \iy$, is a subspace of $L^p(G)$ and for $p=1$ 
it is even a Banach
algebra for the norm $\no{\cdot}_{\o,1}$.
\begin{proposition}\label{restriction}
Let $G$ be a locally compact group and let $S$ be a closed normal
subgroup of $G$. Let $\o$ be a continuous weight on $G$ such that
the inverse of its restriction to $S$ is integrable with respect
to the Haar measure on $S$. Let $f,g: G\to \C$ be two continuous
functions on $G$, such that $\o\cd g$ is uniformly bounded and
such that $ f\in L^1(G,\o)$. Let $h:=f*g\in L^1(G,\o)$. Then for
every $t\in G$, the function $\tilde{h}(t)$ is in $L^1(S)$ and the
mapping $G\times G\to L^1(S); (t,u)\mapsto {}^u{\tilde{h}(t)}$ is
continuous.
\end{proposition}

\begin{proof}
Since $\o$ is a weight, we have that $ \o(s)\leq \o(u)\o(u\inv s)$, 
i.e.~$ \frac 1 {\o(u\inv s)}\leq \frac{\o(u)}{\o(s)},\ \ s,u\in G. $ Hence, for $t\in G,\ s\in S,$
\begin{eqnarray*}
\va{\tilde h(t)(s)}&=& \va{\int_G f(u)g(u\inv ts)\, du} =\va{
\int_G
f(t u)g(u\inv s)\, du}\\
&\leq&\int_G \va {f(tu)}\va{g(u\inv s)}\frac{\o(u\inv s)}{\o(u\inv
s)}\, du \leq\int_G \va {f(tu)}\o(u)\va{g(u\inv s)}\frac{\o(u\inv
s)}{\o( s)}\, du\\
\end{eqnarray*}
 and so
\begin{eqnarray}
\no{\tilde h(t)}_1 &\leq& \int_S\int_G \va {f(tu)}\o(u)\va{g(u\inv
s)}\frac{\o(u\inv s)}{\o( s)}\ du ds\nonumber\\
&\leq& \int_S\int_G \va {f(tu)}\o(u)
\frac{\no{g}_{\o,\iy}}{\o(s)}\ du ds\nonumber
\\
\label{omega}&\leq& \int_S\int_G \o(t\inv)\va
{f(tu)}\o(tu)\frac{\no{g}_{\o,\iy}} {\o (s)}\, du ds=\o(t)\no {f
}_{\o,1}\no {g}_{\o,\iy} \no{(\frac 1 \o)\res S}_1
\end{eqnarray}
So for every $t\in G$, the function $\tilde{h}(t)$ is in $L^1(S)$.
Furthermore, for $t,t'\in G$, by \eqref{omega},
$$\no{h(t)-h(t')}_1\
\leq \int_S\int_G \va {f(tu)-f(t'u)}\o(u)\frac{\no{g}_{\o,\iy}}
{\o (s)}\, du ds$$
$$
\leq \no {(\lambda(t\inv)f-\lambda({t'}\inv)f) }_{\o,1}\no {g}_{\o,\iy} \no{(\frac 1 \o)\res S}_1,
$$ 
where $\lambda$ denotes left translation by elements of $G$. Since left translation in
$L^1(G,\o)$ and conjugation in $L^1(S)$ are continuous, it follows
that the mapping $(t,u)\mapsto {}^u\tilde h(t)$ from $G\times G$
to $L^1(S)$ is continuous too.
\end{proof}

\noindent Let as in \eqref{gauss-estimate} $\delta$ denote the Carath\'eodory
distance associated to our sub-Laplacian $L$ on $G$ and
$(h_t)_{t>0}$ its heat kernel. Then the function
$\o_d(g):=e^{d\delta(x,e)},\ g\in G$, $d\in\R_+$, defines a weight
on $G$. Since we have the Gaussian estimate
$$ 
\vert h_{t}(g)\vert \leq C_t e^{-C_t \d(g,e)^{2}}\textrm{, for all }g\in G, t>0, 
$$
it follows that:
\begin{equation}\label{gauss}
h_t \in L^1(G,\o_d)\cap L^\iy(G,\o_d) \hb{ for every } t>0 \hb{ and }d>0.
\end{equation}
\begin{proposition}\label{decrease}
Let $G$ be the semidirect product of a connected compact Lie group 
$K$ acting on 
an exponential solvable Lie group $S$. Then there exists a constant $d>0$, 
such that
$\frac 1 {\o_d}|_S $ is in $L^1(S)$.
\end{proposition}
\begin{proof}
Let $U$ be a compact symmetric neighborhood of $e$ in $G$
containing $K$. Since $S$ is connected, we know that
$G=\cup_{k\in\N} U^k$. This allows us to define
$\tau_U=\tau:G\to \N$ by:
$$\tau(x)=\min\{k\in\N|\, x\in U^k\}.$$
Then $\tau$ is sub-additive and defines thus a distance on $G$,
which is bounded on compact sets. Since $\tau$ is clearly
connected in the sense of \cite{vsc}, it follows that $\tau$ and
the Carath\'eodory distance $\delta$ are equivalent at
infinity, i.e.
$$1+\tau(x)\leq D(1+ \delta(x))\leq D'(1+\tau(x)),\ x\in G.$$
We choose now a special compact neighborhood of $e$ in the
following way. We take our $K$-invariant scalar-product on $\s$,
the unit-ball $B_\a$ in $\a$ and the unit-ball $B_\n \in \n$. Both
balls are $K$-invariant. Let $U_\a=\exp{B_\a}$ and
$U_{\n}=\exp{B_{\n}}$. Then $U=K U_\a U_\n\cap U_\n U_\a K$ is a
compact symmetric neighborhood of $e$. Let us give a rough
estimate of the radia of the "balls" $U^l,\ l\in\N$. For
simplicity of notation, we shall denote all the positive constants
which will appear in the following arguments (and which will be
assumed to be integers, if necessary) by $C$.\\
Let $k_ia_i n_i\in K {U_\a }{U_\n},\ i=1,\cds ,l$ and
$g:=\Pi_{i=1}^l k_i a_in_i$. We have
$$
g=\Pi_{i=1}^l k_i a_in_i=(\Pi_{i=1}^lk_ia_i)
\left((k_2a_2\cds k_l a_l)\inv n_1 (k_2 a_2 \cds k_l a_l)\cds (k_l a_l)n_{l-1}(k_l a_l) n_l\right).
$$
Since $U_\a$ is $K$-invariant, it follows that:
$$
g=\Pi_{i=1}^l k_i a_is_i= k' a'\Pi_{i=1}^l (a''_i k''_i) n_i (a''_i k''_i)\inv,
$$
where $k', k''_1,\cds k''_l\in K$, 
$a'\in U_\a^l,\ a''_1\in U_\a^{l-1},\cds, a''_{l-1}\in U_a$. 
Hence there exists $X_1,\cds, X_l\in B_\a$, such that 
$$
a'=\exp{X_1}\cds \exp{X_l}=\exp{(X_1+\cds +X_l)}\exp{q_{l}(X_1,\cds,X_l)}
$$ 
for some element $q_{l}(X_1,\cds,X_l)\in \n\cap \s_0$. Since $\s_0$ is a nilpotent
Lie algebra we have that 
$\no{q_l(X_1,\cds,X_l)}\leq C(1+l)^C,\ l\in \N.$ Hence
$$
a'\in \exp{(l B_a)}\exp{[C(1+l)^C B_\n]}\subset \exp{(l B_\a)}U_\n ^{C (1+l)^C}.
$$
Furthermore, because $U_\a$ is compact, $\sup_{a\in U_\a}\opno {\hb{ Ad}(a)}\leq C<\iy$ 
and so $(a''_i k''_i) n_i (a''_i k''_i)\inv \in \exp{ C^{(l-i)} B_\n}\subset U_\n^{C^{l-i}}\ (i=1,\cds,l).$ 
Finally for some integer constants $C$,
\begin{eqnarray}
\nonumber g=k' a'\Pi_{i=1}^l (a''_i k''_i) n_i (a''_i k''_i)\inv
& \in & K \exp{ l U_a} U_\n^{C(1+l)^C} \left(\Pi_{i=1}^{l-1} U_\n^{C^{l-i}}\right) U_\n \\
\nonumber & \subset& K \exp{ l U_a} U_\n^{ C(1+l)^C+\sum_{i=1}^{l-1} C^{l-i}+1} \\
\nonumber & \subset & K\exp{ l U_a} U_\n^{C^l} \\
&\subset & K \exp{ l U_a} \exp{{C^l} B_\n}
\end{eqnarray}
Hence for any $g\in G$, for $\tau_U(g)=l$, we have that $g\in (K U_\a U_\n)^l$ and so, 
denoting by $\hb{Log}: S\to \s$ the inverse map of
$\hb{exp}: \s \to S$, $g=k_ga_gn_g$, with $k_g\in K$,
$a_g\in \exp {\a}$, $\no{\Log{a_g}}$$\leq l=\tau_U(g)$ and $n_g\in N$ with $ \no{\Log{n_g}}$$\leq C^l$, 
i.e.~$\log(1+\no{\Log{n_g}})$$\leq C l=C \tau_G(g)$. Whence for our weight
$\o_d$, ($d\in\R_+$), we have that:
\begin{displaymath}
\begin{array}{c}
\o_d(g)=e^{d \delta(g)}\geq C e^{d C\tau_U(g)}\geq C e^{d C\left(\no {\Log{a_g}}+
\log(1+\no{\Log{n_g}})\right)}\\
= C e^{d C\no {\Log{a_g}}}(1+\no{\Log {n_g}})^{\mathit{dC}}.
\end{array}
\end{displaymath}
Therefore, for $d$ big enough,
$$
\begin{array}{c}
\displaystyle{\int_S \frac{1}{\o_d(s)}ds=\int_{\a}\int_\s\frac{1} {\o_d(\exp X \exp Y)}\, dY dX}\\
\displaystyle{\leq C\int_{\a}\int_\s e^{-d C\no {X}}\frac 1{(1+\no{Y})^{dC}}\, dY dX <\iy.}
\end{array}
$$

\end{proof}

\begin{proposition}\label{continuous-compact}
Let $T$ be a compact topological space and let $k:T\times T\to \ca K(\ca H)$ 
be a continuous mapping into the space of compact
operators on a Hilbert space $\ca H$. Let $\mu$ be a Borel
probability measure on $T$. Then the linear mapping
$$
\begin{array}{c}
\displaystyle{K:L^2(T,\ca H)\to L^2(T,\ca H),}\\
\displaystyle{K\xi(t):=\int_T k(t,u)\xi(u)\, du,\ t\in T,\ \xi\in L^2(T, \ca H),}
\end{array}
$$ 
is compact too.
\end{proposition}

\begin{proof}
We show that $K$ is the norm-limit of a sequence of operators of
finite rank. Let $\varepsilon>0$. Since $T$ is compact and $k$ is
continuous, there exists a finite partition of unity of $T\times T$ 
consisting of continuous non-negative functions
$(\v_i)_{i=1}^N$, such that
$\opno{k(t,t')-k(u,u')}<\frac\e 2$ for every $(t,t'), (u,u')$
contained in the support $\v_i $. Choose for $i=1,\cds, N$ an
element $(t_i,t'_i)$ in supp $\v_i$. Since $k(t_i,t'_i)$ is a
compact operator, we can find a bounded endomorphism $F_i$ of 
$\ca H$ of finite rank, such that $\opno{k(t_i,t'_i)-F_i}<\frac \e 2$,
hence $\opno{k(t,t')-F_i}< \e $ for every $(t,t')\in$ supp $\v_i$,
$i=1,\cds,N$. The finite rank operator $F_i$ has the expression 
$F_i= \sum_{k=1}^{N_i} P_{\eta_{i,k},\eta'_{i,k}}$, where for
$\eta,\eta'\in \ca H$, $P_{\eta,\eta'}$ denotes the rank one
operator $P_{\eta,\eta'}(\eta'')=\sp{\eta''}{\eta'} \eta,\ \eta''\in\ca H$.

We approximate the continuous functions $\v_i$ uniformly on
$T\times T$ up to an error of at most $\frac \e{R}$ by tensors
$\psi_i=\sum_{j=1}^{M_i} \v_{i,j}\otimes\v'_{i,j}\in C(T,\R_+)\otimes C(T,\R_+)$ for some $R>0$ to be
determined later on. Let $K_\e$ be the finite rank operator

$$
K_\e=\sum_{i=1}^{N}\sum_{j=1}^{M_i}
\sum_{k=1}^{N_i}P_{\v_{i,j}\otimes\eta_{i,k},\ \v'_{i,j}\otimes\eta'_{i,k}}.
$$

In order to estimate the difference $ {K-K_\ve}$, we let first
$K_{\ve,1}$ be the kernel operator with kernel
$k_{\ve,1}(s,t)=\sum_{i=1}^N \v_i(s,t)F_i$. Then for $\xi\in L^2(T,\ca H)$
$$\no {K_{\ve,1}\xi-K\xi}_2^2=\int_T\no {\sum_{i=1}^N\int_T \v_i(s,t)(k(s,t)-F_i)\xi(t)dt}^2\,ds$$
$$\leq \int_T(\sum_{i=1}^N\int_T \v_i(s,t)\e \no{\xi(t)}dt)^2\,ds
=\int _T(\int_T\e \no{\xi(t)}dt)^2\,ds\leq \e^2 \no\xi^2,$$
hence $\opno{K-K_{\e,1}}\leq\e.$ Moreover
$$\no{(K_{\e,1}-K_{\e})\xi}^2
=\int_T\no { \int_T\sum_{i=1}^N(\v_i(s,t)-\sum_{j=1}^{M_i}\v_{i,j}(s)\v'_{i,k}(t))F_i\xi(t)\,dt}^2\,ds$$
$$\leq\int_T({\int_T \sum_{i=1}^N\frac \e{R} \opno {F_i}\no {\xi(t)}dt})^2\,ds
\leq \frac{\e^2}{R^2}(\sum_{i=1}^N\opno {F_i})^2 \no\xi^2.$$
So, if we let $R=\frac 1{1+\sum_{i=1}^N\opno {F_i}}$, then
$$\opno{K-K_\e}\leq \opno{K-K_{\e,1}}+\opno{K_{\e,1}-K_\e}\leq 2\e.$$
\end{proof}
Let now $\pi$ be an isometric  representation of the
group $S$ on a Banach space $X$  and denote by 
$\rho^p:={\mathrm ind}^G_{S}\pi$ be  the corresponding induced 
representation of $G$
on $L^p(G,X; \pi).$ Here we follow notation of Section \ref{induced}.\\
Let $h$ be in $L^1(G)$ and assume furthermore that
$\tilde{h}(g)\in L^1(S)$ for all $g\in G$ and that the mapping
$\tilde{h}: G\to L^1(S)$ is continuous. Then the operator
$\rho^p(h)$ is a kernel operator, whose kernel $k(t,u),\, t, u\in G,$ is 
given by:
\begin{equation}\label{kernels}
k(t,u)=\Delta_G(u^{-1})\pi({}^u \tilde h(t u\inv))
\end{equation}
(in the notations of Proposition \ref{restriction}). 
Indeed, for $\xi\in L^p(G,X ;\pi)$, $t\in G$,
\begin{eqnarray*}
[\rho^p(h)\xi](t)&=&\int_G h(g)\xi(g\inv t)\, dg=\int_G
\Delta_G(g\inv)
h(tg\inv)\xi(g)\, dg \\
 &=&\int_{G/S}\int_S \Delta_G(s\inv g\inv)
h(ts\inv g\inv)\xi(g s)\, ds dg \\
 &=& \int_{G/S}\int_S \Delta(g\inv)\Delta_S(s\inv )
h(ts\inv g\inv)\xi(g s)\, ds dg \\
 &=& \int_{G/S}\int_S \Delta(g\inv)
h(tg \inv (g s g\inv))\pi(s)\xi(g )\, ds dg\\
 &=&\int_{G/S} \Delta(g\inv)\pi({ }^g\tilde
h(tg\inv))\xi(g)\,dg.
\end{eqnarray*}
Moreover the kernel $k$ satisfies the following covariance property under $S$:
\begin{equation}\label{cov}
k(ts,us')=\pi(s^{-1})k(t,u)\pi(s'),\qquad t,\, u\in G,\, s,\,s'\in S.
\end{equation}

\begin{proposition}\label{compact}
Let $G$ be the semidirect product of a connected compact Lie group 
$K$ acting on an exponential solvable Lie group $S$. 
Let $(\pi,\ca H)$ be an irreducible unitary representation of the
normal closed subgroup $S$ of $G$ whose Kirillov-orbit
$\O_\pi=\O\subset \s^*$ is closed. Let $\rho=\ind S G \pi.$ Then
the operator $\rho(h_t)$ is compact for every $t>0$.
\end{proposition}
\begin{proof}
By the relation \eqref{gauss} the function $h_t$ is in $L^1(G,\o_d)\cap L^\iy(G,\o_d)$ 
for $t $ and $d>0$. Furthermore we have that $h_t=h_{t/2}*h_{t/2}$. Hence by the Propositions
\ref{decrease} and \ref{restriction} the mapping 
$G\times G\to L^1(S),\, (s,u)\mapsto {}^u\tilde h_t(su\inv)$, is continuous and
so the operator valued kernel function $k(s,u):=\Delta_G(u^{-1})\pi({}^u\tilde{h}_t(su\inv))$ 
is continuous too. It follows from the preceeding discussion
that the
$k$ is just the integral  kernel of the operator $\rho(h_t)$. The fact that
the Kirillov
orbit of
$\pi\in\widehat S$ is closed in $\s^*$ implies that for every
$\v\in L^1(S)$, the operator $\pi(\v)=\int_S f(s) \pi(s)ds $ is
compact (see \cite{leptin-ludwig} and \cite{hlm}). Hence $k(s,u)$
is compact for every $(s,u)\in G\times G$ and in particular for every
$(s,u)\in K\times K$. We apply
Proposition \ref{continuous-compact} to the restriction of $k$ to $K\times K$. The
related kernel operator on $L^2(K,\mathcal{H})$ is then compact. Now, since
$\pi$ is
unitary, the restriction map to $K$ is an isometric isomorphism from
$L^2(G,\mathcal{H};\pi)$ onto
$L^2(K,\mathcal{H}),$ and we thus see that $\rho(h_t)$ is
compact too.
\end{proof}

\subsection{Proof of Theorem \ref{maintheorem3}}\label{proof3}
\setcounter{equation}{0}

We now turn to the proof of Theorem \ref{maintheorem3}, which
follows closely the notation and argumentation in \cite{hlm}. In the
sequel, we always make the following
\begin{quote}
{\bf Assumption.} $\ell\in\s^*$ satisfies Boidol's
condition \eqref{boidolcond}, and $\O(\ell)|_\n$ is closed.\\
Moreover, we assume that $p\in [1,\infty[,\ p\ne 2,$ is
fixed.
\end{quote}
Since $\ell$ satisfies \eqref{boidolcond}, the stabilizer $\s(l)$ is not
contained in $\n$. Let $\nu$ be the real character of $\s$, 
which has been defined in \cite[Section 5]{hlm}, 
trivial on $\n$ and different from $0$ on $\s(\ell)$. We denote by
$\pi_\ell={\rm ind}_P^S
\chi_\ell$ the irreducible unitary representation of $S$ associated to
$\ell$ by the
Krillov map; here $P=P(\ell)$ denotes a suitable polarizing subgroup for
$\ell,$
and $\chi_\ell$ the character $\chi_\ell(p):=e^{i\ell(\log p)}$ of $P.$\\

For any complex number $z$ in the strip
$$\Sigma:=\{\zeta\in\C: |\Im\zeta|<1/2\},$$
let $\De_z$ be the complex character of $S$ given by
$$\De_z(\exp X):=e^{-i z\nu(X)}, \quad X\in\s,$$
and $\chi_z$ the unitary character
$$\chi_z(\exp X):=e^{-i\Re z\nu(X)}, \quad X\in\s.$$
If we define $p(z)\in ]1,\infty[$ by the equation
\begin{equation}
\Im z=1/2 -1/p(z),
\end{equation}
it is shown in \cite{hlm} that the representation $\pi_\ell^z$,
given by
\begin{equation}\label{pitau}
\pi_\ell^z(x):=\De_z(x)\pi_\ell(x)=\chi_z(x)\pi_\ell^{{\ol
{p(z)}}}(x), \quad x\in G,
\end{equation}
is an isometric representation on the mixed $L^{ p}$-space
$L^{{\ol {p(z)}}}(S/P,\ell).$ Here, $\pi_\ell^{{\ol
{p(z)}}}$ denotes the $\ol{p(z)}$-induced representation of $S$ on
$L^{{\ol {p(z)}}}(S/P,\ell)$ defined in \cite{hlm}, where $\ol {p(z)}$ is
a multi-index of the form $(p(z),\dots,p(z),2,\dots,2).$\\
Observe that for $\tau\in\R,$ we have $p(\tau)=2,$ and
$\pi_\ell^\tau= \chi_\tau\otimes\pi_\ell$ is a unitary
representation on $L^{\ol 2}(S/P,\ell)$. Moreover,
\begin{equation}
\pi_\ell^{\tau}\simeq\pi_{\ell-\tau\nu,}
\end{equation}
since the mapping $f\mapsto {\ol \chi}_\tau f$ intertwines the
representations $\chi_\tau \otimes\pi_\ell$ and
$\pi_{\ell-\tau\nu}.$\\

We take now for $z\in\Sigma$ the $p(z)$-induced
representation
$\rho^z_\ell:=\ind {p(z),S} G \pi^z_\ell$ of $G$ which acts on the
space
$$L^{\ol {p(z)}}(G/P,\ell):=L^{p(z)}(G,L^{\ol {p(z)}}(S/P,\ell);\pi^z_\ell).$$
Let us shortly write
$$L^{\ol p}:={L^{\ol {p(z)}}(G/P,\ell)}, \quad 1\le p<\infty,$$
for the space of $\rho^z_\ell$.\\
%
%
\noindent We can extend the character $\De_z,\, z\in\Sigma,$ of
$S$ to a function on $G$ by letting
$$\De_z(k a n):=\De_z(a n)= e^{-i z\nu (\Log a)},\, k\in
K, a\in A, n\in N.$$
Since $\nu$ is trivial on $\n$ and since $kak\inv \in a N$ for all
$k\in K, a\in A$, we have that
$$\De_z( k a n k')=\De_z(a n),\, \, k,k'\in K,a\in A, n\in N,$$
and in particular $\De_z$ is a character of $G$.

\medskip
\noindent Define the operator $T(z),\ z\in \Sigma$, by:
$$T(z):=\rho_\ell^z(h_1).$$
Then by the relations \eqref{kernels} and \eqref{pitau}, for
$z\in\Sigma $ and $\xi\in L^{\ol p}$, (since $\De_z$ is
$K$-invariant)
\begin{eqnarray}
\nonumber T(z)\xi(k) &=& \int_K \pi_\ell^z({}^{k'}\tilde
h_1(k{k'}\inv))\xi(k') dk'\\
\nonumber &=& \int_{K} \pi_\ell((\De_z{\res S}){}^{k'}\tilde
h_1(k{k'}\inv))\xi(k') dk' \\
\nonumber &=&\int_{K} \pi_\ell({}^{k'}
\widetilde{(\De_z h_1)}(k{k'}\inv))\xi(k') dk' \\
\nonumber &=&[\rho_\ell(\De_z h_1)](\xi(k)).
\end{eqnarray}
Hence
\begin{equation}\label{T(z)}
T(z)=\rho^z_\ell(h_1)=\rho_\ell (\De_z h_1),\ z\in \Sigma.
\end{equation}
Since by \eqref{gauss}, for every continuous character $\chi$ of
$G$ which is trivial on $N$ the function $\chi h_1$ is in
$L^1(G)$, it follows from \cite[Corollary 5.2 and Proposition 3.1]{hlm}
that the
operator $T(z)$ leaves $L^{\ol q}$ invariant for every $1\le q<\infty,$ and
is bounded
on all these spaces. Moreover, by Proposition \ref{compact}, $T(\tau)$ is
compact for
$\tau\in\R$. From here on we can proceed exactly as in the proof of
\cite[Theorem 1]{hlm}, provided that we can prove a ``Riemann-Lebesgue"
type lemma like \cite[Theorem 2.2]{hlm} in our
present setting, since $G=K\ltimes S$ is amenable.

We must show that $T(\tau)$ tends to 0 in the operator norm if
$\tau$ tends to $\iy$ in $\R$. The condition we have imposed on
the coadjoint orbit $\O$ of $\ell$, namely that the restriction of
$\O$ to $\n$ is closed, tells us that 
$\lim_{\tau\to\iy} \O+\tau\nu =\iy$ in the orbit space, which means that 
$\lim_{\tau\to\iy}\opno{\pi_{\ell+\tau \nu}(f)} =0$ for every $f\in L^1(S)$. Now, by
\eqref{kernels} the operator $T(\tau)=\rho^\tau_\ell(h_1)$ is a
kernel operator whose kernel $K_\tau$ has values in the bounded
operators on $\ca H_\ell$. the kernel $K_\tau$ is given by:
$$K_\tau(k,k')= \int_S\De_\tau( s )h_1(k\inv s{k'}\inv ) \pi_\ell(s) ds=\pi_\ell ^\tau(h_1(k,k')),$$ 
where $h_1(k,k')$ is the function on $S$ defined
by $h_1(k,k')(s):=h_1(k s {k'}\inv)$. Hence 
$$\lim_{\tau\to
\iy}\opno{\pi_\ell ^\tau (h_1(k,k'))}=0$$ 
for every $k,k'\in K$. Moreover for $k,k'\in K$,
$$\opno{\pi_\ell ^\tau(h_1(k,k'))}\leq \no{ h_1(k,k')}_1\leq \sup_{k''\in K} \no{\tilde{h}_1(k'')}_1 .$$ 
We know from 
Proposition \ref{restriction} that, for every $k''\in K$,
$$\no{h_1(k'')}_1\leq
\no{\o_d|_K}_\iy \no {h_{1/2}}_{\o_d,1}\no {h_{1/2}}_{\o_d,\iy} \no{(\frac 1 \o_d)|_S}_1,$$ 
which is finite by Proposition \ref{decrease} and relation 
\eqref{gauss} (if $d$ is big enough) . Hence, by Lebesgue's dominated convergence theorem,
we see that:
$$\lim_{\tau\to \iy}\int_K\int_K \opno{\pi^\tau_\ell(h_1(k,k')}^2 \, dk d k'=0.$$ 
This shows that:
$$\lim_{\tau\to\iy} \opno{\rho^\tau_\ell(h_1)}=0.$$


\providecommand{\bysame}{\leavevmode\hbox to3em{\hrulefill}\thinspace}
\providecommand{\MR}{\relax\ifhmode\unskip\space\fi MR }
\providecommand{\MRhref}[2]{%
  \href{http://www.ams.org/mathscinet-getitem?mr=#1}{#2}
}
\providecommand{\href}[2]{#2}

\vspace{1cm}

\noindent\emph{Universit\'e de Metz, Math\'ematiques, Ile du
Saulcy, 57045 Metz Cedex, France\\ e-mail:
ludwig@poncelet.univ-metz.fr}

\vspace{0.6cm}

\noindent\emph{Mathematisches Seminar, C.A.-Universit\"at Kiel,
Ludewig-Meyn-Str.4, D-24098 Kiel, Germany\\
e-mail: mueller@math.uni-kiel.de}

\vspace{0.6cm}

\noindent\emph{IRMA - UFR de Math\'ematique et d'Informatique de 
Strasbourg, 7, rue Ren\'e Descartes, 67084 Strasbourg Cedex, France\\
e-mail: souaifi@math.u-strasbg.fr} 

\end{document}